\documentclass[twocolumn]{autart}    

\usepackage{graphicx}          
\usepackage[utf8]{inputenc}
\usepackage{textcomp}

\usepackage{amsmath}

\usepackage{siunitx}
\usepackage{longtable,tabularx}
\usepackage{subfigure}
\usepackage{algorithm}
\usepackage{algpseudocode}
\usepackage{xcolor}
\usepackage{comment}
\usepackage{upgreek}
\usepackage{graphicx}
\usepackage{booktabs}
\usepackage{hyperref}
\usepackage{setspace}
\usepackage{mathtools}
\usepackage{amsfonts}
\usepackage{natbib}
\usepackage{amssymb}
\usepackage{bigints}
\usepackage{xspace}
\usepackage{multirow}
\usepackage{array}



\graphicspath{{images/}}

\DeclareMathOperator*{\argmax}{\arg\!\max}

\begin{document}

\begin{frontmatter}

\title{Fixed-Time State Transfer via Pontryagin Extremals\thanksref{footnoteinfo}}

\thanks[footnoteinfo]{This paper was not presented at any conference.}

\author[KAIST]{Juho Bae}\ead{juhobae@mit.edu}, 
\author[KAIST]{Han--Lim Choi\thanksref{corrautinfo}}\ead{hanlimc@kaist.ac.kr}

\address[KAIST]{Korea Advanced Institute of Science and Technology, Daejeon, 34141, Republic of Korea}

\thanks[corrautinfo]{Corresponding author.}

\begin{keyword}                           
Boundary value problem; Pontryagin maximum principle; Reachability; Dubins car; Motion planning                
\end{keyword}                             

\begin{abstract}                         
	This paper concerns the problem of fixed-time transition between two states of nonlinear systems (i.e., the point-to-point steering problem). We propose a formulation applicable to a broad class of nonlinear systems, and show that it is theoretically complete in the sense that it admits a solution if and only if the target state is reachable from the initial state. When the target state is reachable, we prove that a solution can always be constructed by concatenation of two Pontryagin extremals, one generated by the original dynamics $f$ from the initial state and the other generated by the inverted dynamics $-f$ from the terminal state. This allows the problem to be formulated as a two-point boundary value problem (TPBVP) of extremals, where the solution existence to the formulated TPBVP is equivalent to that of the original problem. The theoretical developments are applied to curves with prescribed curvature bounds in $\mathbb{R}^3$, thereby extending the recent works on Dubins car to dimension three. We prove that to construct a curvature-bounded path in $\mathbb{R}^3$ with prescribed length and boundary conditions, it suffices to consider the trajectories that are concatenations of CSC, CCC, their subsegments, and H, where C denotes a circular arc with maximum curvature, S a straight line segment, and H a certain class of helicoidal arcs with constant curvature. Numerical demonstrations are conducted on a nonlinear dynamics example, and on curvature-bounded paths in $\mathbb{R}^2$ and $\mathbb{R}^3$.
\end{abstract}

\end{frontmatter}

%
\section{Introduction} \label{sec:01}
This paper concerns the task of transition between two points of a nonlinear system within a fixed time interval. More precisely, consider the process in $\mathbb{R}^{n_{\boldsymbol{\chi}}}$
\begin{equation} \label{eq:dynamics}
	\dot{\boldsymbol{\chi}} = f(\boldsymbol{\chi}, \boldsymbol{u}), \quad \boldsymbol{\chi}(0) = \boldsymbol{\chi}_i
\end{equation}
where $\boldsymbol{\chi} \in \mathbb{R}^{n_{\boldsymbol{\chi}}}$, $\boldsymbol{u} \in \mathbb{R}^{n_{\boldsymbol{u}}}$, and $f$ is $C^1$ in $\mathbb{R}^{n_{\boldsymbol{\chi}} + n_{\boldsymbol{u}}}$. The initial value is represented by $\boldsymbol{\chi}_i \in \mathbb{R}^{n_{\boldsymbol{\chi}}}$. Let $\mathcal{F}$ be the family of all measurable control inputs $\boldsymbol{u}(t)$ on $[0, T]$ where each $\boldsymbol{u}(t)$ lies in $\Omega$, a compact restraint set in $\mathbb{R}^{n_{\boldsymbol{u}}}$. The point-to-point steering problem is formulated as follows.

\textbf{Problem P1}. \textit{Given $\boldsymbol{\chi}_i, \boldsymbol{\chi}_f \in \mathbb{R}^{n_{\boldsymbol{\chi}}}$, find a control $\boldsymbol{u} \in \mathcal{F}$ and the corresponding trajectory $\boldsymbol{\chi}(t)$ such that $\boldsymbol{\chi}(0) = \boldsymbol{\chi}_i$ and $\boldsymbol{\chi}(T) = \boldsymbol{\chi}_f$.}

Solving problem \textbf{P1} in differentially flat systems is an easy task through algebraic means without explicit integration~(\cite{fliess1995}). However, it remains challenging when the system is not differentially flat. To tackle such dynamics, a two-point boundary value problem (TPBVP) formulation of problem \textbf{P1} was introduced in ~\cite{Graichen2005} and~\cite{Graichen2008} based on partial feedback linearization. It was proposed to introduce free parameters to provide additional degree of freedom to make the boundary value problem (BVP) of the internal dynamics solvable. \cite{Kvitko2017,kvitko2019one,kvitko2019method,kvitko2020} treat point-to-point steering problem of control affine systems as a sequence of stabilization problems. \cite{Fetisov2020} reduces the problem into linear controllable subsystems through partial feedback linearization and provides sufficient condition for solution existence. However, the following fundamental concerns still remain in the formulations of these works.
\begin{enumerate}
	\item Existence of solution to problem \textbf{P1}, 
	\item Given the solution existence to problem \textbf{P1}, is the formulation guaranteed to have a solution?
\end{enumerate}
The first question, which relates to the concept of \textit{reachability}, has a rich history~(\cite{bansal2017}). Despite their computational cost, theoretical frameworks~(\cite{lygeros2004}) allow for determining solution existence of problem \textbf{P1} through Hamilton-Jacobi (HJ) reachability analysis. Nevertheless, the second question remains unsolved. The formulations proposed in the aforementioned studies either lack guarantee of the existence of solution even when the target state $\boldsymbol{x}_f$ is reachable~(\cite{Graichen2005,Graichen2008}), or require overly restrictive assumptions on the system dynamics~(\cite{Kvitko2017,kvitko2019one,kvitko2019method,kvitko2020}), or provide only sufficient conditions stronger than reachability for solution existence~(\cite{Fetisov2020}).

The papers cited above do not specify a certain dynamics, but another major stream of studies on point-to-point steering problem exists in \textit{Dubins car} context. Dubins car dynamics represent curves in $\mathbb{R}^2$ with a prescribed curvature bound~(\cite{dubins1957}), and serves as a good approximation to the dynamics of many practical vehicles~(\cite{Ding2019}). Thus, the problem of finding a feasible path of a Dubins car with a prescribed length and boundary conditions (i.e., prescribed departure and arrival locations and directions) rises in multiple practical applications. Although this problem is essentially problem \textbf{P1} with dynamics specified as Dubins car, numerous research were conducted utilizing the geometric property of curvature bound.

While there exists multiple approaches (e.g., curve fitting, homotopy) in literature (\cite{xu1999,yao2017}), a common strategy is to first find a curve of minimum length satisfying the boundary conditions using the conclusions of~\cite{dubins1957} and appropriately elongate it to the desired length. Some examples include~\cite{Schumacher2003} and~\cite{Meyer2015}. However, the aforementioned works did not address the two previously discussed concerns, first, the existence of a trajectory of given length and boundary conditions, and second, whether the desired trajectory can always be found by the proposed elongation strategies when existence is assumed. It is only recently that complete solution addressing these issues was proposed in~\cite{chen2023} and~\cite{Rao2024}. The authors present explicit conditions for solution existence leveraging the geometric description of the reachable set of Dubins car in~\cite{patsko2003} and~\cite{patsko2022}. While the contributions of these two works are significant, the arguments made are hard to directly handle different types of boundary conditions (e.g., when certain state variable is unconstrained at the boundary, one has to redo the proofs as in~\cite{Ding2019} and~\cite{chen2023}), or to be generalized to higher dimensions or different dynamics as they heavily leveraged the geometric properties of Dubins car.

The limitation of such geometry-driven elongation approaches is particularly evident in three dimensions. Indeed, planning of curvature-bounded paths in $\mathbb{R}^3$ is a key primitive in robotic motion planning, particularly for fixed-wing aerial vehicles with constrained turning radius~(\cite{hota2010}). However, while the planar problem was extensively studied for decades before being completely resolved by the recent geometric constructions~(\cite{chen2023,Rao2024}), attempting a similar geometric elongation approach in $\mathbb{R}^3$ is profoundly difficult. Although time-optimal trajectories are characterized in~\cite{sussmann1995}, the elongation strategies that are effective in the planar setting do not admit a straightforward extension to $\mathbb{R}^3$, primarily due to the intricate coupling of torsion and curvature inherent in helicoidal curves. This highlights the necessity of a fundamentally different formulation that does not rely on geometric elongation techniques.

To this end, we propose a formulation for the point-to-point steering problem via Pontryagin extremals with an existence guarantee: whenever problem \textbf{P1} is feasible, the proposed formulation admits a solution. The contributions of this paper are summarized as follows: First, we establish the fundamental connection between problem \textbf{P1} and certain set-valued maps that contain the boundary of reachable sets. This connection implies our main result; \textit{given the existence of solution to problem \textbf{P1}, there always exists a solution consisting of concatenation of two Pontryagin extremals; one generated by the original dynamics $f$ from the initial state and the other generated by the inverted dynamics $-f$ from the terminal state.} Thus, problem \textbf{P1} can be formulated as a TPBVP in a manner that is free from the second concern; existence of a solution to this TPBVP is equivalent to that of the original point-to-point steering problem. Moreover, the proposed approach is capable of addressing different boundary conditions when certain states are unconstrained. Second, we apply the theoretical developments to extend the existing results on Dubins car~(\cite{chen2023,Rao2024}) to curvature-bounded paths in $\mathbb{R}^3$. Since our theoretical developments bypass the need for geometric elongation, the arguments in planar setting extends in a unified manner to spatial setting. Specifically, we prove that for construction of such paths under prescribed length and boundary conditions, it suffices to consider the trajectories that are concatenations of CSC, CCC, their subsegments, and H, where C denotes a circular arc with maximum curvature, S a straight line segment, and H a certain class of smooth helicoidal arcs with constant curvature. This facilitates the use of root-finding algorithms (e.g., the Levenberg-Marqardt algorithm) to find the desired path, even in the presence of unconstrained boundary states as well.
%
\section{Setup and Preliminaries} \label{sec:02}
\subsection{Notations and conventions}
Let us commence by introducing the notations frequently used throughout the paper. We use $\|\cdot\|$ to represent Euclidean norm of vectors and spectral norm of matrices. The set of all interior points of a set $S$ is denoted by $int(S)$ and the set of all boundary points by $bd(S)$. Set complement is represented as $S^\complement$. Distance between a point $x$ and a set $S$ is defined and denoted as $d(x, S) = \inf\limits_{y \in S} | x - y |$ where distance between points are usual Euclidean distance. In the space of nonempty compact subsets of $\mathbb{R}^n$, the distance between two sets $S_1$ and $S_2$ is defined by Hausdorff distance as $d_H(S_1, S_2) = \max \left\{ \sup\limits_{x_1 \in S_1} d(x_1, S_2), \sup\limits_{x_2 \in S_2} d(x_2, S_1) \right\}$. Continuity of nonempty compact set--valued maps are with respect to Hausdorff metric induced topology unless explicitly stated otherwise. The symbol \textit{a.e.} $t$ stands for \textit{almost all $t$} condition.
\subsection{Preliminaries on Pontryagin maximum principle and reachability}
Throughout this paper, we introduce the following assumptions on $f$ in Eq.~\eqref{eq:dynamics}. 
\begin{assum} \label{assum:01}
	\hfill
	\begin{itemize}
		\item There exists a uniform bound $\|\boldsymbol{\chi}(t)\| < M$ for all responses of the dynamical system Eq.~\eqref{eq:dynamics} on $[0, T]$.
		\item For each $\boldsymbol{\chi}$, the set $V(\boldsymbol{\chi}) = \{ f(\boldsymbol{\chi}, \boldsymbol{u}) : \boldsymbol{u} \in \Omega \}$ is convex.
	\end{itemize}
\end{assum}
Here, a \textit{uniform bound} means that the choice of the upper bound $M$ is independent of the choice of the controller $\boldsymbol{u} \in \mathcal{F}$ and time. Compactness of $V(\boldsymbol{\chi})$ follows trivially by compactness of $\Omega$ and continuity of $f$. These assumptions along with the others that will be introduced throughout this paper are met by useful classes of systems, particularly control affine systems under an appropriate convex restraint set $\Omega$. Nonetheless, we maintain a general form of the dynamics to facilitate future expandability. Under these assumptions, the \textit{reachable set} $\mathcal{G}(f, t, \boldsymbol{\chi}_i)$ of $f$ at $t \in [0, T]$ is defined as follows.
\begin{defn}
\begin{equation}
	\mathcal{G}(f, t, \boldsymbol{\chi}_i) \equiv \{ \boldsymbol{\chi}(t) : \boldsymbol{u} \in \mathcal{F} \}
	\nonumber
\end{equation}
\end{defn}
The notion of reachability can be related to problem \textbf{P1} as follows: finding the reachable set $\mathcal{G}(f, T, \boldsymbol{\chi}_i)$ is to determine the set of all $\boldsymbol{\chi}_f \in \mathbb{R}^{n_{\boldsymbol{\chi}}}$ such that problem \textbf{P1} has a solution.

Under Assumption~\ref{assum:01}, Filippov's compactness theorem for reachable sets, applied with the standard cut-off argument for a priori bounded responses, implies that $\mathcal{G}(f,t,\boldsymbol{\chi}_i)$  is compact~(\cite{AgrachevSachkov2004}). Moreover, $\mathcal{G}(f, t, \boldsymbol{\chi}_i)$ varies continuously in $t$~(\cite{Lee1967}). The following theorem provides a necessary condition on the trajectories reaching the boundary of $\mathcal{G}(f, t, \boldsymbol{\chi}_i)$.
\begin{thm}{(\cite{Lee1967})} \label{thm:PMP}
	Suppose $\boldsymbol{u}^*(t) \in \mathcal{F}$ have a response $\boldsymbol{\chi}^*(t)$ such that $\boldsymbol{\chi}^*(t_f) \in bd(\mathcal{G}(f, t_f, \boldsymbol{\chi}_i))$ for some $t_f \in [0, T]$. Then there exists a nontrivial adjoint system $\boldsymbol{p}^*(t)$ on $[0, t_f]$ such that 
	\begin{equation} \label{eq:thm1-1}
		\dot{\boldsymbol{p}}^*(t) = -\boldsymbol{p}^*(t)^T \frac{\partial f}{\partial \boldsymbol{\chi}}( \boldsymbol{\chi}^*(t), \boldsymbol{u}^*(t) ),
	\end{equation}
	\begin{equation} \label{eq:thm1-2}
	\begin{split}
		& \mathcal{H}( \boldsymbol{p}^*(t), \boldsymbol{\chi}^*(t), \boldsymbol{u}^*(t) ) = \max\limits_{\boldsymbol{u} \in \Omega} \mathcal{H}( \boldsymbol{p}^*(t), \boldsymbol{\chi}^*(t), \boldsymbol{u} ) \\
		& \text{for a.e. } t, 
	\end{split}
	\end{equation}
	where $\mathcal{H}( \boldsymbol{\chi}, \boldsymbol{p}, \boldsymbol{u} ) \equiv \boldsymbol{p}^T f( \boldsymbol{\chi}, \boldsymbol{u} )$.
\end{thm}
We define the set of all solution triplets $\{ \boldsymbol{\chi}^*, \boldsymbol{p}^*, \boldsymbol{u}^* \}$ of Eqs.~\eqref{eq:thm1-1} and \eqref{eq:thm1-2} as follows along with the set of endpoints of such solutions.
\begin{defn}
Given dynamics $f$, $t_f \in [0, T]$, and $\boldsymbol{\chi}_i \in \mathbb{R}^{n_{\boldsymbol{\chi}}}$, we define the following sets.
\begin{itemize}
	\item $S(f, t_f, \boldsymbol{\chi}_i) \equiv \{ \{ \boldsymbol{\chi}^*, \boldsymbol{p}^*, \boldsymbol{u}^* \} : \text{Nontrivial solution of } \\ \text{Eqs.~\eqref{eq:thm1-1} and \eqref{eq:thm1-2} on } [0, t_f] \}$
	\item $E(f, t_f, \boldsymbol{\chi}_i) \equiv \{ \boldsymbol{\chi}^*(t_f) : \{ \boldsymbol{\chi}^*, \boldsymbol{p}^*, \boldsymbol{u}^* \} \in S(f, t_f, \boldsymbol{\chi}_i) \}$
\end{itemize}
\end{defn}
Additionally, we introduce the following terminologies. While~\cite{Lee1967} uses the term \textit{maximal} instead of extremal, we insist on this terminology to avoid confusion with other literatures.
\begin{defn}
Any trajectory $\boldsymbol{\chi}$ of $f$ on $[0, t_f]$ is called \textit{extremal} if it belongs to $\left\{ \boldsymbol{\chi} : \{ \boldsymbol{\chi}, \boldsymbol{p}, \boldsymbol{u} \} \in S(f, t_f, \boldsymbol{\chi}_i) \right\}$. An extremal trajectory $\boldsymbol{\chi}$ on $[0, t_f]$ is \textit{extendable} to $T$ if there exists $\{ \boldsymbol{\chi}', \boldsymbol{p}', \boldsymbol{u}' \} \in S(f, T, \boldsymbol{\chi}_i)$ such that $\boldsymbol{\chi}(t) = \boldsymbol{\chi}'(t)$ for all $t \in [0, t_f]$.
\end{defn}
Since $\mathcal{G}(f, t_f, \boldsymbol{\chi}_i)$ is compact under Assumption~\ref{assum:01}, all of its boundary points are in $\mathcal{G}(f, t_f, \boldsymbol{\chi}_i)$ and are therefore reachable at time $t_f$. Consequently, the following corollary holds.
\begin{cor} \label{cor:01}
	$bd(\mathcal{G}(f, t_f, \boldsymbol{\chi}_i)) \subseteq E(f, t_f, \boldsymbol{\chi}_i)$
\end{cor}
It should be noted that identifying $E(f, t_f, \boldsymbol{\chi}_i)$ does not completely determine the boundary of $\mathcal{G}(f, t_f, \boldsymbol{\chi}_i)$ in general, because Theorem~\ref{thm:PMP} only serves as a necessary condition but not sufficient. Instead, the following statement allows to exclude some points in $E(f, t_f, \boldsymbol{\chi}_i)$ that lie interior to $\mathcal{G}(f, t_f, \boldsymbol{\chi}_i)$; the points in $E(f, t_f, \boldsymbol{\chi}_i)$ that can be reached by an alternative trajectory that is not extremal (i.e., violates the condition Eq.~\eqref{eq:thm1-2}) lie interior to $\mathcal{G}(f, t_f, \boldsymbol{\chi}_i)$. This is because Theorem~\ref{thm:PMP} states there exists a nontrivial solution of Eqs.~\eqref{eq:thm1-1} and~\eqref{eq:thm1-2} for `every' trajectory reaching the boundary of $\mathcal{G}(f, t_f, \boldsymbol{\chi}_i)$. This argument was leveraged during the reachability analysis in~\cite{patsko2003} and~\cite{bae2025reachability}.
\subsection{Related works on curves in $\mathbb{R}^2$ and $\mathbb{R}^3$ with prescribed curvature bounds}
Shortest paths with a prescribed curvature bound have been extensively studied in the literature. The planar case was first addressed in~\cite{dubins1957}, concluding that the minimizers have structures of CSC, CCC, or their respective subsegments. Analogous conclusion was reproduced later by applying PMP to the following Dubins car dynamics~(\cite{Pecsvaradi1972}).
\begin{subequations} \label{eq:dubins_dynamics}
	\begin{align}
		\boldsymbol{\chi}_{2D} &:= \begin{bmatrix}x \\ y \\ \gamma\end{bmatrix}, 
		\quad 
		f_{2D}(\boldsymbol{\chi}_{2D},u_{2D}) \;:=\;
     	\begin{bmatrix}\cos\gamma\\ \sin\gamma\\ 0\end{bmatrix}
     	+ \begin{bmatrix}0\\ 0\\ 1\end{bmatrix}u_{2D} \\[2pt]
		\dot{\boldsymbol{\chi}}_{2D} &= f_{2D}(\boldsymbol{\chi}_{2D},u_{2D}), \quad u_{2D} \in [-1, 1]
	\end{align}
\end{subequations}
The state variables $x$ and $y$ represent the geometric coordinates in $\mathbb{R}^2$ and $\gamma \in \mathbb{R}$ represents the heading angle. Without loss of generality, we assume unit speed, unit curvature bound, and zero initial state denoted as $\boldsymbol{\chi}_{2D, i} := \boldsymbol{\chi}_{2D}(0) = (0, 0, 0)$. The desired terminal state is represented as $\boldsymbol{\chi}_{2D, f} \in \mathbb{R}^3$. Since the signed curvature equals to $u$ wherever defined, the control restraint set is defined as $\Omega = [-1, 1]$. The major conclusion on reachable sets of Dubins car that will be utilized throughout this paper is as follows. We note that~\cite{patsko2003} provides additional necessary conditions, but we only leverage the following result for simplicity.
\begin{prop} \label{prop:Dubins2D}
	(\cite{patsko2003}) Any boundary point of $\mathcal{G}(f_{2D}, t, \boldsymbol{\chi}_{2D, i})$ can be reached by the curves belonging to one of the following classes: CSC, CCC, and their subsegments.
\end{prop}
The proof is done by first constructing $E(f_{2D}, t, \boldsymbol{\chi}_{2D, i})$ through solving Eqs.~\eqref{eq:thm1-1} and~\eqref{eq:thm1-2}, and then identifying some family of extremals such that their endpoints lie interior to the reachable set. Since certain extremals are neglected during this process, and that we omit some necessary conditions of extremals in Proposition~\ref{prop:Dubins2D}, it is important to note that there is no inclusion relationship between the set of CSC and CCC curves and $S(f_{2D}, t, \boldsymbol{\chi}_{2D, i})$.

The optimal control theory of curvature-bounded paths in $\mathbb{R}^3$ was studied in~\cite{sussmann1995} by introducing the following dynamics 
\begin{subequations} \label{eq:3d_dynamics}
	\begin{align}
		\boldsymbol{\chi}_{3D} &:= \begin{bmatrix}\boldsymbol{x} \\ \boldsymbol{y} \end{bmatrix}, 
		\quad 
		f_{3D}(\boldsymbol{\chi}_{3D},\boldsymbol{u}_{3D}) \;:=\;
     	\begin{bmatrix}\boldsymbol{y}\\ \boldsymbol{y} \times \boldsymbol{u}_{3D} \end{bmatrix} \\[2pt]
		\dot{\boldsymbol{\chi}}_{3D} &= f_{3D}(\boldsymbol{\chi}_{3D},\boldsymbol{u}_{3D}), \quad \boldsymbol{u}_{3D} \in \mathbb{B}^3
	\end{align}
\end{subequations}
where $\boldsymbol{x} \in \mathbb{R}^3$ represents the geometric coordinates, $\boldsymbol{y} \in S^2$ the velocity, and $\mathbb{B}^3$ the closed unit ball in $\mathbb{R}^3$. Similar to Dubins car, we assume unit speed and curvature bound, and denote the initial state as $\boldsymbol{\chi}_{3D, i} := \boldsymbol{\chi}_{3D}(0) = (0, 0, 0, 1, 0, 0)$. The desired terminal state is represented as $\boldsymbol{\chi}_{3D, f} \in \mathbb{R}^3 \times S^2$. \cite{sussmann1995} provides a detailed description of the structure of shortest paths between two states of the above dynamics through application of PMP. To introduce the conclusions of~\cite{sussmann1995}, let us first consider the smooth helicoidal arcs that have constant curvature of $1$ and nonvanishing torsion $\tau$ satisfying the following ODE in Eq.~\eqref{eq:torsionODE} for some constant $\zeta \in \mathbb{R}$.
\begin{equation} \label{eq:torsionODE}
	\ddot{\tau} = \frac{3\dot{\tau}^{2}}{2\tau} - 2\tau^3 + 2\tau - \zeta \tau \sqrt{|\tau|}
\end{equation}
We denote the class of these helicoidal arcs by H. Conclusion of~\cite{sussmann1995} was that the optimal solutions to the minimum time problem are of CSC, CCC, their respective subsegments, or H with $\zeta \geq 0$. 

By utilizing the developments of~\cite{sussmann1995},~\cite{bae2025reachability} addressed the reachability of curvature-bounded paths in $\mathbb{R}^3$ through similar steps in~\cite{patsko2003}. The conclusion that will be utilized throughout this paper is summarized in Proposition~\ref{prop:Dubins3D}. As in Proposition~\ref{prop:Dubins2D}, we only leverage the argument in Proposition~\ref{prop:Dubins3D} while the authors of~\cite{bae2025reachability} provide additional necessary conditions for extremals. Similarly, we note that there is no inclusion relationship between the set of trajectories outlined in Proposition~\ref{prop:Dubins3D} and $S(f_{3D}, t, \boldsymbol{\chi}_{3D, i})$.
\begin{prop} \label{prop:Dubins3D}
	(\cite{bae2025reachability}) Any boundary point of $\mathcal{G}(f_{3D}, t, \boldsymbol{\chi}_{3D, i})$ can be reached by the curves belonging to one of the following classes: CSC, CCC, their subsegments, and H.
\end{prop}
%
\section{Main Results} \label{sec:03}
Building upon the previous works outlined, this section presents the central arguments of this paper. We first introduce the notion of \textit{en route reachability} to present the connections between PMP and point-to-point steering problems. The main conclusions are outlined in Theorems~\ref{thm:main} and~\ref{thm:extremal}. The subsection following presents the TPBVP formulation of point-to-point steering problem when the structure of extremals are determined by the value of initial costate vector.
\subsection{Connections between reachability and point-to-point steering problem}
We introduce the concept of \textit{en route reachable tube}, which represents the set of all states that can be reached during the transition from initial to the terminal state.
\begin{defn}
	(En route reachable tube)
	Given endpoint $\boldsymbol{\chi}_f \in \mathbb{R}^{n_{\boldsymbol{\chi}}}$, let $\mathbb{X}(f, T, \boldsymbol{\chi}_i, \boldsymbol{\chi}_f)$ denote the set of all state solutions of problem \textbf{P1}. En route reachable tube is the set of all points that the trajectories in $\mathbb{X}(f, T, \boldsymbol{\chi}_i, \boldsymbol{\chi}_f)$ can pass through, defined as follows.
	\begin{equation}
		\mathcal{B}(f, T, \boldsymbol{\chi}_i, \boldsymbol{\chi}_f) := \bigcup\limits_{\boldsymbol{\chi} \in \mathbb{X}(f, T, \boldsymbol{\chi}_i, \boldsymbol{\chi}_f)} \boldsymbol{\chi}\left([0, T]\right)
		\nonumber
	\end{equation}
\end{defn}
It is easy to see that $\mathcal{B}(f, T, \boldsymbol{\chi}_i, \boldsymbol{\chi}_f)$ is nonempty if and only if problem \textbf{P1} has a solution. Subsequently, the set of en route reachable points at time $t$ is defined as follows.
\begin{defn}
	(En route reachable set)
	For $t \in [0, T]$, the set of all points in $\mathbb{R}^{n_{\boldsymbol{\chi}}}$ such that the trajectories between $\boldsymbol{\chi}_i$ and $\boldsymbol{\chi}_f$ can pass through at time $t$ is defined as follows.
	\begin{equation}
		\mathcal{R}(t) \equiv \left\{ \boldsymbol{\chi}(t): \boldsymbol{\chi} \in \mathbb{X}(f, T, \boldsymbol{\chi}_i, \boldsymbol{\chi}_f) \right\} 
		\nonumber
	\end{equation}
\end{defn}
Recall that $\mathcal{G}(f, t, \boldsymbol{\chi}_i)$ is the set of all reachable points at $t$, through the trajectories of dynamics $f$, initiating from $\boldsymbol{\chi}_i$. We denote the set of all such trajectories on $[0, T]$ as $\mathbb{X}_1$ so that $\mathcal{G}(f, t, \boldsymbol{\chi}_i) = \{ \boldsymbol{\chi}_1(t) : \boldsymbol{\chi}_1 \in \mathbb{X}_1 \}$. In other words, each $\boldsymbol{\chi}_1 \in \mathbb{X}_1$ satisfies $\boldsymbol{\chi}_1(t) = \boldsymbol{\chi}_i + \int_0^t f(\boldsymbol{\chi}_1(\tau), \boldsymbol{u}(\tau)) d\tau$ for some $\boldsymbol{u} \in \mathcal{F}$. Similarly, $\mathcal{G}(-f, t, \boldsymbol{\chi}_f)$ is the set of all reachable points at $t$, through the trajectories of dynamics $-f$, initiating from $\boldsymbol{\chi}_f$. We define $\mathbb{X}_2$ as the set of all such trajectories on $[0, T]$ so that $\mathcal{G}(-f, t, \boldsymbol{\chi}_f) = \{ \boldsymbol{\chi}_2(t) : \boldsymbol{\chi}_2 \in \mathbb{X}_2 \}$ and $\boldsymbol{\chi}_2(t) = \boldsymbol{\chi}_f + \int_0^t -f(\boldsymbol{\chi}_2(\tau), \boldsymbol{u}(\tau)) d\tau$.

From these definitions, it is easy to see that any $\boldsymbol{\chi} \in \mathbb{X}(f, T, \boldsymbol{\chi}_i, \boldsymbol{\chi}_f)$ is a concatenation of some $\boldsymbol{\chi}_1 \in \mathbb{X}_1$ and $\boldsymbol{\chi}_2 \in \mathbb{X}_2$, in the sense that $\boldsymbol{\chi}_1(\tau) = \boldsymbol{\chi}_2(T-\tau)$ for some $\tau \in [0, T]$. Moreover, $\mathcal{R}(t)$ is the set of all such `concatenation points' at time $t$, so that $\mathcal{R}(t) = \{ \boldsymbol{\chi}_1(t) : \boldsymbol{\chi}_1(t) = \boldsymbol{\chi}_2(T-t), t \in [0, T], \boldsymbol{\chi}_1 \in \mathbb{X}_1, \boldsymbol{\chi}_2 \in \mathbb{X}_2 \}$. It then follows $\mathcal{R}(t) = \mathcal{G}(f, t, \boldsymbol{\chi}_i) \cap \mathcal{G}(-f, T-t, \boldsymbol{\chi}_f)$ and $\bigcup\limits_{t \in [0, T]} \mathcal{R}(t) = \mathcal{B}(f, T, \boldsymbol{\chi}_i, \boldsymbol{\chi}_f)$. The former implies that $\mathcal{R}(t)$ is compact. A useful yet nontrivial property of $\mathcal{R}(t)$ is continuity, which further implies compactness of $\mathcal{B}(f, T, \boldsymbol{\chi}_i, \boldsymbol{\chi}_f)$.
\begin{lem} \label{lemma:continuity}
	$\mathcal{R}(t)$ is continuous in $t$ on $[0, T]$.
\end{lem}
\begin{pf*}{Proof.}
	If we pick a point $\boldsymbol{p}_1 \in \mathcal{R}(t_1)$ for some $t_1 \in [0, T]$, there exists a control $\boldsymbol{u} \in \mathcal{F}$ and the corresponding trajectory $\boldsymbol{\chi}$ such that $\boldsymbol{\chi}(t_1) = p_1$, $\boldsymbol{\chi}(0) = \boldsymbol{\chi}_i$, and $\boldsymbol{\chi}(T) = \boldsymbol{\chi}_f$. From the uniform bound assertion in Assumption~\ref{assum:01} and compactness of $\Omega$, there exists a uniform bound: $\left|f(\boldsymbol{\chi}, \boldsymbol{u})\right| < m$. Then for any $\varepsilon > 0$ and $t_2 \in [0, T]$ such that $|t_1 - t_2| < \frac{\varepsilon}{m}$, $|\boldsymbol{\chi}(t_1) - \boldsymbol{\chi}(t_2)| < m |t_1 - t_2| < \varepsilon$. Thus, each point $\boldsymbol{p}_1 = \boldsymbol{\chi}(t_1) \in \mathcal{R}(t_1)$ lies within $\varepsilon$--distance from some point $\boldsymbol{\chi}(t_2) \in \mathcal{R}(t_2)$. But analogous arguments can be applied to any $\boldsymbol{p}_2 \in \mathcal{R}(t_2)$ to prove that it is within $\varepsilon$--distance from some point in $\mathcal{R}(t_1)$. Thus, $d_H(\mathcal{R}(t_1), \mathcal{R}(t_2)) < \varepsilon$ if $|t_1 - t_2| < \frac{\varepsilon}{m}$. Hence, the map $t \mapsto \mathcal{R}(t)$ is a continuous map from $[0, T]$ to the metric space of nonempty compact sets in $\mathbb{R}^{n_{\boldsymbol{\chi}}}$. \qed
\end{pf*}
Before presenting our main theorem, we introduce the following key topological lemma. The example followed by articulates the related insights.
\begin{lem} \label{lemma:bubble}
	Suppose $S(t)$ and $A(t)$ are nonempty, compact-valued, and continuous set-valued maps on an interval $[a, b]$. For any continuous function $\boldsymbol{x}(t): [a, b] \mapsto \mathbb{R}^{n_{\boldsymbol{\chi}}}$, suppose $\boldsymbol{x}(a) \in int\left( S(a) \right)$ and $\boldsymbol{x}(b) \notin int\left( S(b) \right)$. If $bd\left( S(t) \right) \subseteq A(t)$ for all $t \in [a, b]$, then there exists $t' \in [a, b]$ such that $\boldsymbol{x}(t') \in A(t')$.
\end{lem}
\begin{pf*}{Proof.}
	Assume the converse that for all $t \in [a, b]$, $\boldsymbol{x}(t)$ never touches $A(t)$ and thus $bd\left( S(t) \right)$. We define two sets $T_1 \equiv \left\{ t \in [a, b] : \boldsymbol{x}(t) \in int\left( S(t) \right) \right\}$ and $T_2 \equiv \left\{ t \in [a, b] : \boldsymbol{x}(t) \in S(t)^\complement \right\}$. From the hypothesis, it follows $T_1 \cap T_2 = \emptyset$ and $T_1 \cup T_2 = [a, b]$.
	
	Let us denote $d\left( \boldsymbol{x}(t), S(t) \right)$ by $d(t)$ on $[a, b]$. Then continuity of $\boldsymbol{x}(t)$ and $S(t)$ implies continuity of $d(t)$. Since $S(t)$ is compact, $\boldsymbol{x}(t) \in S(t)^\complement$ if and only if $d\left( t \right) > 0$. Moreover, $\boldsymbol{x}(t) \notin bd\left(S(t)\right)$ implies that $\boldsymbol{x}(t) \in int\left(S(t)\right)$ if and only if $d\left( t \right) = 0$. Consequently, $T_1 = d^{-1}\left( \{0\} \right)$ and continuity of $d$ implies that $T_1$ is closed in $[a, b]$.
	
	Now, suppose $T_1$ has a boundary point $\tilde{t}$. Then there exists a sequence $t_n \rightarrow \tilde{t}$ such that $\boldsymbol{x}(t_n) \in S(t_n)^\complement$. Moreover, since $T_1$ is closed, $\tilde{t} \in T_1$ and $d\left( \tilde{t} \right) = 0$. Then from the relationship $bd\left( S(t_n) \right) \subseteq A(t_n)$ and continuity of $A(t)$ and $d(t)$, it follows 
	\begin{equation}
		\begin{split}
			0 &\leq d\left(\boldsymbol{x}(\tilde{t}), A(\tilde{t})\right) \\
			& = \lim\limits_{n\rightarrow \infty} d\left(\boldsymbol{x}(t_n), A(t_n)\right) \\
			& \leq \lim\limits_{n\rightarrow \infty} d\left(\boldsymbol{x}(t_n), bd\left(S(t_n)\right)\right) \\
			& = \lim\limits_{n\rightarrow \infty} d\left(\boldsymbol{x}(t_n), S(t_n)\right) \\
			& = \lim\limits_{n\rightarrow \infty} d(t_n) \\
			& = d(\tilde{t}) \\
			& = 0.
		\end{split}
	\end{equation}
	Consequently, $d(\boldsymbol{x}(\tilde{t}), A(\tilde{t})) = 0$ and compactness of $A(\tilde{t})$ implies $\boldsymbol{x}(\tilde{t}) \in A(\tilde{t})$.
		
	If the closed set $T_1$ has no boundary point, then it is either empty or $[a, b]$. However, $T_1$ is nonempty because $\boldsymbol{x}(a) \in int\left(S(a)\right)$. Moreover, $b \notin T_1$ because $\boldsymbol{x}(b) \notin int\left( S(b) \right)$. Thus, $T_1$ must have a boundary point. This completes the proof. \qed
\end{pf*}
\begin{exmp}
	Consider the region $S(t) = \{\boldsymbol{x} \in \mathbb{R}^2 : \|\boldsymbol{x}\| \leq 1 - t \}$ enclosed by a `bubble' $A(t) := \partial S(t) = \{\boldsymbol{x} \in \mathbb{R}^2 : \|\boldsymbol{x}\| = 1 - t \}$, where $t \in [0, 1]$. The bubble shrinks to a point at time $t = 1$. If a particle $\boldsymbol{p}(t) \in \mathbb{R}^2$ moves continuously starting at $\boldsymbol{p}(0) \in int\left( S(0) \right)$, then Lemma~\ref{lemma:bubble} states that there exists a time $t' \in [0, 1]$ such that $\boldsymbol{p}(t') \in A(t')$. In other words, the particle $\boldsymbol{p}(t)$ eventually touches the bubble $A(t)$.
\end{exmp}
Building upon the previous lemmas, our main theorem is stated as follows. This theorem highlights the intrinsic connection between reachability and point-to-point steering problem.
\begin{thm} \label{thm:main}
	Suppose $A_1(t)$ and $A_2(t)$ are nonempty, compact-valued, and continuous set--valued maps defined on $[0, T]$. For all $\boldsymbol{\chi}_1 \in \mathbb{X}_1$ and $\boldsymbol{\chi}_2 \in \mathbb{X}_2$, suppose the followings hold.
	\begin{itemize}
		\item $bd\left(\mathcal{G}(f, t, \boldsymbol{\chi}_i)\right) \subseteq A_1(t) \subseteq \mathcal{G}(f, t, \boldsymbol{\chi}_i)$ and \\ $bd\left(\mathcal{G}(-f, T-t, \boldsymbol{\chi}_f)\right) \subseteq A_2(t) \subseteq \mathcal{G}(-f, T-t, \boldsymbol{\chi}_f)$ for all $t \in [0, T]$.
		\item If $\boldsymbol{\chi}_1(t') \in A_1(t')$ $\left( \text{resp. } \boldsymbol{\chi}_2(t') \in A_2(T - t') \right)$ for some $t' \in [0, T]$, there exists a control $\boldsymbol{u}_1$ $\left( \text{resp. } \boldsymbol{u}_2 \right)$ on $[t', T]$ such that its response satisfy $\boldsymbol{\chi}_1(t) \in A_1(t)$ $\left( \text{resp. } \boldsymbol{\chi}_2(t) \in A_2(T - t) \right)$ for all $t \in [t', T]$.
	\end{itemize}
	Then for any $\boldsymbol{x} \in \mathcal{B}(f, T, \boldsymbol{\chi}_i, \boldsymbol{\chi}_f)$, there exists $\tau \in [0, T]$, $\boldsymbol{\chi}_1 \in \mathbb{X}_1$, and $\boldsymbol{\chi}_2 \in \mathbb{X}_2$ such that 
	\begin{enumerate}
		\item $\boldsymbol{\chi}_1(\tau) = \boldsymbol{\chi}_2(T-\tau) \in A_1(\tau) \cap A_2(\tau)$,
		\item $\boldsymbol{x}$ is in $\boldsymbol{\chi}_1\left([0,\tau]\right)$ or $\boldsymbol{\chi}_2\left([0,T-\tau]\right)$.
	\end{enumerate}
\end{thm}
\begin{pf*}{Proof.}
	The idea is as follows. For $\boldsymbol{x}\in\mathcal{B}(f,T,\boldsymbol{\chi}_i,\boldsymbol{\chi}_f)$, consider the ``bubble'' $A(t)=A_1(t)\cup A_2(t)$ which includes the set of boundary points of $\mathcal{R}(t)$. Since $\mathcal{R}(t)$ is continuous in $t$, we use Lemma~\ref{lemma:bubble} to identify the earliest time at which $\boldsymbol{x}$ hits $A(t)$. W.l.o.g.\ this contact occurs on $A_1$, and the invariance assumption allows to extend the corresponding trajectory while staying on $A_1$. Applying Lemma~\ref{lemma:bubble} once more with $A_2$ yields the existence of time $\tau$ where the two extremals meet, giving the desired concatenation.

	Suppose $\boldsymbol{x} \in \mathcal{R}(\overline{t})$ for some $\overline{t} \in [0, T]$. Define a set-valued map $A(t) = A_1(t) \cup A_2(t)$. From the assumption, $A(t)$ is nonempty, compact-valued, and continuous. It is evident from the relationship $\mathcal{R}(t) = \mathcal{G}(f, t, \boldsymbol{\chi}_i) \cap \mathcal{G}(-f, T-t, \boldsymbol{\chi}_f)$ that $bd\left(\mathcal{R}(t)\right) \subseteq bd\left( \mathcal{G}(f, t, \boldsymbol{\chi}_i) \right) \cup bd\left( \mathcal{G}(-f, T-t, \boldsymbol{\chi}_f) \right)$, which implies that $bd\left(\mathcal{R}(t)\right) \subseteq A(t)$ for all $t \in [0, T]$.

	First, suppose $\boldsymbol{x} \in bd\left( \mathcal{R}(\overline{t}) \right)$ so that $\boldsymbol{x}$ lies on the boundary of $\mathcal{G}(f, \overline{t}, \boldsymbol{\chi}_i)$ or $\mathcal{G}(-f, T-\overline{t}, \boldsymbol{\chi}_f)$. This implies that $\boldsymbol{x} \in A_1\left(\overline{t}\right)$ or $\boldsymbol{x} \in A_2\left(\overline{t}\right)$. If $\boldsymbol{x}$ lies on the boundaries of both $\mathcal{G}(f, \overline{t}, \boldsymbol{\chi}_i)$ and $\mathcal{G}(-f, T-\overline{t}, \boldsymbol{\chi}_f)$, then $\boldsymbol{x}$ is in both $A_1(\overline{t})$ and $A_2(\overline{t})$, and the claim of the theorem follows trivially. Therefore, we assume without loss of generality that $\boldsymbol{x} \in bd\left(\mathcal{G}(f, \overline{t}, \boldsymbol{\chi}_i)\right)$ and $\boldsymbol{x} \in int\left(\mathcal{G}(-f, T-\overline{t}, \boldsymbol{\chi}_f)\right)$. It trivially follows that $\boldsymbol{x} \in A_1(\overline{t})$.
	
	Second, suppose $\boldsymbol{x} \in int\left( \mathcal{R}(\overline{t}) \right)$. Then $\overline{t} \in (0, T)$ because $\mathcal{R}(0) = \{\boldsymbol{\chi}_i\}$ and $\mathcal{R}(T) = \{\boldsymbol{\chi}_f\}$ are singleton sets with non interior. Then Lemma~\ref{lemma:bubble} applied to $\mathcal{R}(t)$, $A(t)$, and $\boldsymbol{x}$ on $[\overline{t}, T]$ implies the existence of $t'$ such that $\boldsymbol{x} \in A(t')$. Hence, the set $\{ t' : \boldsymbol{x} \in A(t') \}$ is nonempty. Moreover, this set is closed in $[0, T]$\textemdash{}and hence in $\mathbb{R}$\textemdash{}since compactness of $A(t)$ implies that it is the inverse image of the closed set $\{0\}$ of the continuous function $d(\boldsymbol{x}, A(t))$. Consequently, we can consider the smallest $t'$ in $[\overline{t}, T]$ such that $\boldsymbol{x} \in A(t')$, denoted as $t_m$. If there exists $t'' \in (\overline{t}, t_m)$ such that $\boldsymbol{x} \notin int\left( \mathcal{R}(t'') \right)$, Lemma~\ref{lemma:bubble} implies the existence of $t_m' < t_m$ such that $\boldsymbol{x} \in A(t_m')$. This contradicts that $t_m$ is minimum. Thus, $\boldsymbol{x} \in int\left( \mathcal{R}(t) \right)$ on $[\overline{t}, t_m)$. Then $\lim\limits_{t \rightarrow t_m^-} d\left( \boldsymbol{x}, \mathcal{R}(t) \right) = d\left( \boldsymbol{x}, \mathcal{R}(t_m) \right) = 0$ and compactness of $\mathcal{R}(t_m)$ implies that $\boldsymbol{x} \in \mathcal{R}(t_m)$. Thus, we have proved that $\boldsymbol{x}$ is in $A(t_m) = A_1(t_m) \cup A_2(t_m)$ and $\mathcal{R}(t_m)$. If $\boldsymbol{x} \in A_1(t_m) \cap A_2(t_m)$ then the claim of the theorem follows trivially. Thus, we assume without loss of generality that $\boldsymbol{x} \in A_1(t_m)$ and $\boldsymbol{x} \notin A_2(t_m)$. This implies that $\boldsymbol{x} \notin bd\left(\mathcal{G}(-f, T-t_m, \boldsymbol{\chi}_f)\right)$ because $bd\left(\mathcal{G}(-f, T-t_m, \boldsymbol{\chi}_f)\right) \subseteq A_2(t_m)$. But $\boldsymbol{x} \in \mathcal{R}(t_m) = \mathcal{G}(f, t_m, \boldsymbol{\chi}_i) \cap \mathcal{G}(-f, T-t_m, \boldsymbol{\chi}_f)$ implies that $\boldsymbol{x}$ is in the compact set $\mathcal{G}(-f, T-t_m, \boldsymbol{\chi}_f)$, which implies $\boldsymbol{x} \in int\left( \mathcal{G}(-f, T-t_m, \boldsymbol{\chi}_f) \right)$.
	
	Now, the two cases can be handled together as follows. For a trajectory $\boldsymbol{\chi} \in \mathbb{X}(f, T, \boldsymbol{\chi}_i, \boldsymbol{\chi}_f)$ that reaches $\boldsymbol{x}$ at time $t_m$, throughout the previous arguments, we assumed that $\boldsymbol{\chi}(t_m) \in A_1(t_m)$ and $\boldsymbol{\chi}(t_m) \in int\left( \mathcal{G}(-f, T-t_m, \boldsymbol{\chi}_f) \right)$. Let us define a trajectory $\boldsymbol{\chi}_1$ as $\boldsymbol{\chi}_1(t) := \boldsymbol{\chi}(t)$ on $[0, t_m]$. Then $\boldsymbol{\chi}_1(t_m) \in A_1(t_m)$, and one can extend $\boldsymbol{\chi}_1(t)$ on $[t_m, T]$ so that $\boldsymbol{\chi}_1(t) \in A_1(t)$ for all $t \in [t_m, T]$. Since $\boldsymbol{\chi}_1(t)$ is continuous, Lemma~\ref{lemma:bubble} applied to $\mathcal{G}(-f, T-t, \boldsymbol{\chi}_f)$, $A_2(t)$, and $\boldsymbol{\chi}_1(t)$ implies the existence of $\tau$ such that $\boldsymbol{\chi}_1(\tau) \in A_2(\tau)$. Since $A_2(\tau) \subseteq \mathcal{G}(-f, T-\tau, \boldsymbol{\chi}_f)$, existence of $\boldsymbol{\chi}_2 \in \mathbb{X}_2$ such that $\boldsymbol{\chi}_1(\tau) = \boldsymbol{\chi}_2(T-\tau) \in A_1(\tau) \cap A_2(\tau)$ follows from the definition of reachable set. Moreover, $\boldsymbol{x} = \boldsymbol{\chi}(t_m) \in \boldsymbol{\chi}_1\left([0,\tau]\right)$. This completes the proof. \qed
\end{pf*}
A trivial choice of $A_1(t)$ and $A_2(t)$ satisfying the assumptions of Theorem~\ref{thm:main} are the reachable sets themselves, $\mathcal{G}(f, t, \boldsymbol{\chi}_i)$ and $\mathcal{G}(-f, T-t, \boldsymbol{\chi}_f)$, respectively. Then the statement of Theorem~\ref{thm:main} reduces to the definition of $\mathcal{B}(f, T, \boldsymbol{\chi}_i, \boldsymbol{\chi}_f)$: \textit{any point $\boldsymbol{x} \in \mathcal{B}(f, T, \boldsymbol{\chi}_i, \boldsymbol{\chi}_f)$ can be reached by a trajectory that is a concatenation of $\boldsymbol{\chi}_1 \in \mathbb{X}_1$ and $\boldsymbol{\chi}_2 \in \mathbb{X}_2$, where $\boldsymbol{\chi}_1(\tau) = \boldsymbol{\chi}_2(T-\tau)$. (i.e., one reaching $\boldsymbol{x}$ from $\boldsymbol{\chi}_i$ in time $\tau$ with dynamics $f$, and the other reaching $\boldsymbol{x}$ from $\boldsymbol{\chi}_f$ in time $T - \tau$ with dynamics $-f$.)} Thus, this choice does not provide any useful information. Throughout the subsequent discussions, we claim the existence of nontrivial $A_1(t)$ and $A_2(t)$ through Theorem~\ref{thm:PMP}. We commence from the following lemma on the extendability of extremal trajectories and compactness of $E(f, t_f, \boldsymbol{\chi}_i)$. Since the proof involves differential inclusion theory~(\cite{Filippov1988,Aubin1984}), which is not the primary focus of this paper, we place the proof in Appendix~\ref{appendix:A}.
\begin{lem} \label{lemma:extendability}
	Let $U(\boldsymbol{\chi}, \boldsymbol{p}) = \argmax\limits_{\boldsymbol{u} \in \Omega} \mathcal{H}(\boldsymbol{\chi}, \boldsymbol{p}, \boldsymbol{u})$ and define $F(\boldsymbol{\chi}, \boldsymbol{p}) \equiv \begin{Bmatrix} \begin{bmatrix} f(\boldsymbol{\chi}, \boldsymbol{u}) \\ -\nabla_{\boldsymbol{\chi}} f( \boldsymbol{\chi}, \boldsymbol{u} )^T\boldsymbol{p} \end{bmatrix}: \boldsymbol{u} \in U(\boldsymbol{\chi}, \boldsymbol{p}) \end{Bmatrix}$. If $F(\boldsymbol{\chi}, \boldsymbol{p})$ is convex for each $(\boldsymbol{\chi}, \boldsymbol{p})$, then the followings hold for all $t_f \in [0, T]$.
	\begin{itemize}
		\item Any extremal trajectory on $[0, t_f]$ is extendable.
		\item Any infinite sequence of extremal trajectories $\{\boldsymbol{x}_i(t)\}$ on $[0, t_f]$ has a uniformly converging subsequence such that its limit is also an extremal. Consequently, $E(f, t_f, \boldsymbol{\chi}_i)$ is compact.
	\end{itemize}
\end{lem}
\begin{pf*}{Proof.}	
	See Appendix~\ref{appendix:A}. \qed 
\end{pf*}
As stated in Section~\ref{sec:02}, the assumption of Lemma~\ref{lemma:extendability} is also readily met in a large amount of control affine systems with appropriately defined restraint set $\Omega$. Lemma~\ref{lemma:extendability} facilitates the application of Corollary~\ref{cor:01} to find a nontrivial pair of $A_1$ and $A_2$ in the following theorem.
\begin{thm} \label{thm:extremal}
	Suppose the assumption of Lemma~\ref{lemma:extendability} hold. Then $\boldsymbol{x} \in \mathcal{B}(f, T, \boldsymbol{\chi}_i, \boldsymbol{\chi}_f)$ if and only if $\boldsymbol{x}$ can be reached by a concatenation two extremal trajectories: $\boldsymbol{\chi}_1 \in S(f, \tau, \boldsymbol{\chi}_i)$ and $\boldsymbol{\chi}_2 \in S(-f, T - \tau, \boldsymbol{\chi}_f)$ for some $\tau \in [0, T]$.
\end{thm}
\begin{pf*}{Proof.}
	Corollary~\ref{cor:01} and Lemma~\ref{lemma:extendability} imply that the choice of $A_1(t) = E(f, t, \boldsymbol{\chi}_i)$ and $A_2(t) = E(-f, T - t, \boldsymbol{\chi}_f)$ satisfies the assumptions of Theorem~\ref{thm:main}. Then for any $\boldsymbol{x} \in \mathcal{B}(f, T, \boldsymbol{\chi}_i, \boldsymbol{\chi}_f)$, there exist $\tau \in [0, T]$, $\boldsymbol{\chi}_1 \in S(f, \tau, \boldsymbol{\chi}_i)$, and $\boldsymbol{\chi}_2 \in S(-f, T - \tau, \boldsymbol{\chi}_f)$ such that $\boldsymbol{\chi}_1(\tau) = \boldsymbol{\chi}_2(T-\tau)$ and $\boldsymbol{x} \in \boldsymbol{\chi}_1\left([0,\tau]\right) \cup \boldsymbol{\chi}_2\left([0,T-\tau]\right)$. Consequently, concatenation of such $\boldsymbol{\chi}_1$ and $\boldsymbol{\chi}_2$ is in $\mathbb{X}(f, T, \boldsymbol{\chi}_i, \boldsymbol{\chi}_f)$ and reaches $\boldsymbol{x}$. Converse is trivial from definition. \qed
\end{pf*}
We note that in Theorem~\ref{thm:extremal}, $\boldsymbol{\chi}_1(\tau) = \boldsymbol{x}$ does not necessarily hold, i.e., the concatenation point does not have to be $\boldsymbol{x}$ in general. Under Assumption~\ref{assum:01} and the assumption of Lemma~\ref{lemma:extendability}, Theorem~\ref{thm:extremal} states that if the desired terminal state is reachable, then a solution to the associated point-to-point steering problem can always be constructed by concatenating two extremal trajectories, one generated by the original dynamics $f$ from $\boldsymbol{\chi}_i$ and the other generated by the inverted dynamics $-f$ from $\boldsymbol{\chi}_f$. Here, the segment generated by the inverted dynamics can be viewed as an extremal of $f$ traced backward in time from $\boldsymbol{\chi}_f$. Moreover, collecting the images of all such concatenations becomes the en route reachable tube, $\mathcal{B}(f, T, \boldsymbol{\chi}_i, \boldsymbol{\chi}_f)$.
\subsection{TPBVP formulation of point-to-point steering problem}
In this subsection, we present the TPBVP formulation of point-to-point steering problem by leveraging Theorem~\ref{thm:extremal}, particularly to the systems such that the structure of extremals are determined by the initial costate vector. We note that this does not hold for the dynamics of curvature-bounded paths in $\mathbb{R}^2$ and $\mathbb{R}^3$ in Eqs.~\eqref{eq:dubins_dynamics} and~\eqref{eq:3d_dynamics}, while the geometric structures of their extremals are well known~(\cite{patsko2003,sussmann1995}). Thus, we postpone the discussion of the results related to Eqs.~\eqref{eq:dubins_dynamics} and~\eqref{eq:3d_dynamics} to next section.

If $U(\boldsymbol{\chi}, \boldsymbol{p})$ is a singleton set (i.e., one can write the maximizer of Hamiltonian as a function of $\boldsymbol{\chi}$ and $\boldsymbol{p}$) for almost all $t$, then the assumption of Lemma~\ref{lemma:extendability} is satisfied and the extremals are completely determined by the initial value of the costate vector. While extremals can still be uniquely determined by the initial costate value under the generic assumptions outlined in~\cite{Chitour2006} and~\cite{Chitour2008}, even in the presence of singular trajectories~(\cite{Silva2010}), we impose this assumption as it holds for a fairly wide range of practical scenarios~(\cite{taheri2018}). Subsequently, Theorem~\ref{thm:extremal} facilitates the following TPBVP formulation of problem \textbf{P1}. Here, the state and costate variables are reparametrized on $[0, 1]$ and represented by the same notations.

\textbf{Problem P2}. \textit{Let $\boldsymbol{u}_1(\boldsymbol{\chi}_1, \boldsymbol{p}_1) = \argmax\limits_{\boldsymbol{u} \in \Omega} \boldsymbol{p}_1^T f( \boldsymbol{\chi}_1, \boldsymbol{u} )$ and $\boldsymbol{u}_2(\boldsymbol{\chi}_2, \boldsymbol{p}_2) = \argmax\limits_{\boldsymbol{u} \in \Omega} -\boldsymbol{p}_2^T f( \boldsymbol{\chi}_2, \boldsymbol{u} )$. Consider the ODE 
\begin{equation}
	\begin{bmatrix}
		\dot{\boldsymbol{\chi}}_1 \\ \dot{\boldsymbol{p}}_1 \\ \dot{\boldsymbol{\chi}}_2 \\ \dot{\boldsymbol{p}}_2
	\end{bmatrix} = T \begin{bmatrix}
		\tau f(\boldsymbol{\chi}_1, \boldsymbol{u}_1(\boldsymbol{\chi}_1, \boldsymbol{p}_1)) \\ -\tau \nabla_{\boldsymbol{\chi}} f( \boldsymbol{\chi}_1, \boldsymbol{u}_1 )^T\boldsymbol{p}_1 \\ -(1-\tau) f(\boldsymbol{\chi}_2, \boldsymbol{u}_2(\boldsymbol{\chi}_2, \boldsymbol{p}_2)) \\ (1-\tau) \nabla_{\boldsymbol{\chi}} f( \boldsymbol{\chi}_2, \boldsymbol{u}_2 )^T\boldsymbol{p}_2
	\end{bmatrix}
\end{equation}
on interval $[0, 1]$ for an unknown parameter $\tau \in [0, 1]$. Determine $\boldsymbol{p}_1(0)$, $\boldsymbol{p}_2(0)$, and $\tau \in [0, 1]$ such that the following boundary conditions are satisfied.
\begin{equation} \label{eq:bd_conditions}
	\boldsymbol{\chi}_1(0) = \boldsymbol{\chi}_i, \quad \boldsymbol{\chi}_2(0) = \boldsymbol{\chi}_f, \quad \boldsymbol{\chi}_1(1) = \boldsymbol{\chi}_2(1)
\end{equation}
}
Any solution of problem \textbf{P2} is a solution of problem \textbf{P1} while the converse is false in general. From Theorem~\ref{thm:extremal}, the following statements are equivalent.
\begin{enumerate}
	\item $\mathcal{B}(f, T, \boldsymbol{\chi}_i, \boldsymbol{\chi}_f)$ is nonempty, 
	\item Problem \textbf{P1} has a solution, 
	\item Problem \textbf{P2} has a solution.
\end{enumerate}
Thus, the formulation of problem \textbf{P2} is guaranteed to have a solution if the original problem \textbf{P1} does, and the second concern outlined in Section~\ref{sec:01} is not present. Subsequently, off-the-shelf BVP solvers (e.g.,~\cite{rackauckas2017}) can be applied to solve the underdetermined boundary value problem \textbf{P2}. Additioanlly, collecting the images of all state solutions of problem \textbf{P2} becomes the en route reachable tube $\mathcal{B}(f, T, \boldsymbol{\chi}_i, \boldsymbol{\chi}_f)$.

Nevertheless, we note that the fundamental question on the existence of a solution to problem \textbf{P1} still remains. This is because while the existence of solutions to the two problems are equivalent, testing the existence of solution to problem \textbf{P2} still requires computationally expensive methods such as interval analysis~(\cite{Lin2008}) that suffer from curse of dimensionality. General methodology to test the existence of solution to problem \textbf{P1} (i.e., reachability analysis) involves Hamilton-Jacobi partial differential equation~(\cite{lygeros2004}), which has been a notorious computational challenge over the couple of decades. Thus, related works on point-to-point steering problems assume that the final state is reachable~(\cite{Graichen2005,Graichen2008}), or provide only sufficient conditions on solution existence~(\cite{Fetisov2020}), or impose severe restrictions on the dynamics and boundary values~(\cite{Kvitko2017}). Accordingly, we do not tackle reachability analysis but leave this avenue for future work along the discussions in current and the following section on curvature-bounded paths.

The following remark highlights additional flexibility of the proposed formulation. In the subsequent section, we apply the arguments in this section to curves in $\mathbb{R}^3$ with a prescribed curvature bound, thereby extending the results on Dubins car presented by several authors to 3D.
\begin{rem} \label{rem:unconstrained}
	Problem \textbf{P1} with certain variables unconstrained can be analogously addressed through problem \textbf{P2}: if certain boundary states are unconstrained, one can introduce those variables as unknown parameters of the TPBVP. As geometry-driven approaches require to completely redo the proofs to address such variants of point-to-point steering problems (e.g.,~\cite{Ding2019} and~\cite{chen2023}), this carries substantial advantage in simplicity. In principle, one may also consider the case where the terminal time $T$ is unconstrained by treating $T$ as an additional unknown in problem \textbf{P2}. However we do not pursue this extension here, as the present paper focuses on fixed terminal time steering and on the construction guarantee under that setting; allowing free $T$ is left for future work.
\end{rem}
%
\section{Application to curves in $\mathbb{R}^2$ and $\mathbb{R}^3$ with prescribed curvature bounds} \label{sec:04}
In this section, we leverage Theorem~\ref{thm:main} to obtain geometric conclusions on point-to-point steering problem of curvature-bounded paths in $\mathbb{R}^2$ and $\mathbb{R}^3$. For sake of completeness, we note that Assumption~\ref{assum:01} is satisfied for the dynamics in Eq.~\eqref{eq:dubins_dynamics} and Eq.~\eqref{eq:3d_dynamics}. While the choice of $A_1$ and $A_2$ in the previous section was the endpoint mappings of extremals, a more refined choice of $A_1$ and $A_2$ is utilized based on the reachability results in Propositions~\ref{prop:Dubins2D} and~\ref{prop:Dubins3D}.

While the works~\cite{patsko2003} and~\cite{bae2025reachability} on reachability of curvature-bounded paths are based on the application of PMP, the extremals posses singular arcs and the TPBVP formulation in the previous section cannot be applied directly. Moreover, the conclusions in Propositions~\ref{prop:Dubins2D} and~\ref{prop:Dubins3D} are obtained by neglecting some family of extremals such that their endpoints lie interior to the reachable sets\textemdash{}as discussed below Corollary~\ref{cor:01}. Thus, the analyses we present in this section are not covered in the previous discussions.

We note that singular arcs in both $\mathbb{R}^2$ and $\mathbb{R}^3$ cases correspond to the S component~(\cite{patsko2003,sussmann1995}), and the geometric structures of extremals are easily studied thereby. With Proposition~\ref{prop:Dubins2D}, the following theorem can be stated on Dubins car with a straightforward proof.
\begin{thm} \label{thm:dubins2d}
	$\boldsymbol{x} \in \mathcal{B}(f_{2D}, T, \boldsymbol{\chi}_{2D, i}, \boldsymbol{\chi}_{2D, f})$ if and only if $\boldsymbol{x}$ can be reached by a curve of length $T$, and is a concatenation of two trajectories of CSC, CCC, or their subsegments.
\end{thm}
\begin{pf*}{Proof.}
	We first note that the convexity condition of $V(\boldsymbol{\chi})$ in Assumption~\ref{assum:01} is satisfied for the dynamics of Dubins car in Eq.~\eqref{eq:dubins_dynamics}, as it is a control affine system with $\Omega = [-1, 1]$. While the assumption of Lemma~\ref{lemma:extendability} holds as well, our proof only relies on Theorem~\ref{thm:main}. Motivated by Proposition~\ref{prop:Dubins2D}, we define $A_1(t)$ (resp. $A_2(t)$) to be the set of all endpoints of CSC and CCC trajectories (with length of a certain component possibly being zero) of overall length $t$ (resp. $T - t$), initiating from $\boldsymbol{\chi}_{2D, i}$ (resp. $\boldsymbol{\chi}_{2D, f}$). From Proposition~\ref{prop:Dubins2D}, if follows that $bd(\mathcal{G}(f_{2D}, t, \boldsymbol{\chi}_{2D, i})) \subseteq A_1(t) \subseteq \mathcal{G}(f_{2D}, t, \boldsymbol{\chi}_{2D, i})$ and $bd(\mathcal{G}(-f_{2D}, T-t, \boldsymbol{\chi}_{2D, f})) \subseteq A_2(t) \subseteq \mathcal{G}(-f_{2D}, T-t, \boldsymbol{\chi}_{2D, f})$ for each $t$. Extendability is easy to see since elongation of any C or S component at the end is possible through control inputs of $\pm 1$ or $0$. This maintains the elongated trajectory in classes CSC or CCC. Additionally, $A_1(t)$ and $A_2(t)$ vary continuously in $t$ because the endpoints of all CSC and CCC trajectories vary continuously as $t$ varies. It then remains to show that $A_1(t)$ and $A_2(t)$ are compact.

	The endpoint is uniquely determined once the class of the curve and the lengths of the first two components are determined. For instance, if it is given that the curve is of CSC, the winding directions of the C components, and the lengths of the first C and S components, then the endpoint of the CSC curve is determined uniquely. It is easy to see that the relationship between the two lengths and the endpoint is continuous. Moreover, the region of feasible lengths of the first two components is compact, as it consists of nonnegative pair of real numbers with sum $\leq t$ for $A_1(t)$ and sum $\leq T - t$ for $A_2(t)$. Hence, $A_1(t)$ and $A_2(t)$ are continuous images of compact sets, and are therefore compact. Consequently, such choice of $A_1(t)$ and $A_2(t)$ satisfies the conditions outlined in Theorem~\ref{thm:main}. This completes the proof. \qed
\end{pf*}
Theorem~\ref{thm:dubins2d} characterizes the class of curves, consist of maximum six components, required to construct a trajectory of Dubins car of given length and boundary conditions. The corresponding trajectory can be found by solving a nonlinear equation consist of lengths of each components. It is noteworthy that the recent developments in~\cite{chen2023} and~\cite{Rao2024} provide more refined results on Dubins car with lesser number of components through geometric arguments. Nevertheless, our approach is easily extended to $\mathbb{R}^3$, and we focus on the structural similarity of the proofs utilizing Theorem~\ref{thm:main}.
\begin{thm} \label{thm:dubins3d}
	$\boldsymbol{x} \in \mathcal{B}(f_{3D}, T, \boldsymbol{\chi}_{3D, i}, \boldsymbol{\chi}_{3D, f})$ if and only if $\boldsymbol{x}$ can be reached by a curve of length $T$, and is a concatenation of two trajectories of CSC, CCC, their subsegments, or H.
\end{thm}
\begin{pf*}{Proof.}
	The convexity condition of $V(\boldsymbol{\chi})$ in Assumption~\ref{assum:01} is satisfied for $f_{3D}$ in Eq.~\eqref{eq:3d_dynamics}, as it is a control affine system with a convex control restraint set $\Omega = \mathbb{B}^3$. Following the notations of~\cite{sussmann1995}, let us introduce the costate variable $(\boldsymbol{\lambda}, \boldsymbol{\mu}) \in \mathbb{R}^3 \times \mathbb{R}^3$ so that the Hamiltonian of Eq.~\eqref{eq:3d_dynamics} can be written as $\mathcal{H}(\boldsymbol{x}, \boldsymbol{y}, \boldsymbol{\lambda}, \boldsymbol{\mu}, \boldsymbol{u}_{3D}) = \langle \boldsymbol{\lambda}, \boldsymbol{y} \rangle + \langle \boldsymbol{\mu}, \boldsymbol{y} \times \boldsymbol{u}_{3D} \rangle = \langle \boldsymbol{\lambda}, \boldsymbol{y} \rangle + \langle \boldsymbol{u}_{3D}, \boldsymbol{\mu} \times \boldsymbol{y} \rangle$. Then the set of maximizers is $\left\{ \frac{\boldsymbol{\mu} \times \boldsymbol{y}}{|\boldsymbol{\mu} \times \boldsymbol{y}|} \right\}$ if $\boldsymbol{\mu} \times \boldsymbol{y} \neq \boldsymbol{0}$, and $\mathbb{B}^3$ if $\boldsymbol{\mu} \times \boldsymbol{y} = \boldsymbol{0}$. Then the set $F(\boldsymbol{\chi}, \boldsymbol{p})$ defined in Lemma~\ref{lemma:extendability} is convex for each $(\boldsymbol{\chi}, \boldsymbol{p})$ because it is an affine image of the convex set of maximizers (which is either a singleton or $\mathbb{B}^3$). This satisfies the assumption of Lemma~\ref{lemma:extendability}. It was proved in~\cite{sussmann1995} and~\cite{bae2025reachability} that the extremals of Eq.~\eqref{eq:3d_dynamics}\textemdash{}$S(f_{3D}, t, \boldsymbol{\chi}_{3D, i})$\textemdash{}are either H, or certain type of concatenations of C and S components.
	
	Proposition~\ref{prop:Dubins3D} motivates us to set $A_1(t)$ (resp. $A_2(t)$) as the set of all endpoints of CSC, CCC, and H trajectories of overall length $t$ (resp. $T - t$), possibly with length of a certain component being zero, initiating from $\boldsymbol{\chi}_{3D, i}$ (resp. $\boldsymbol{\chi}_{3D, f}$). Now, it remains to show that such choice of $A_1(t)$ and $A_2(t)$ satisfies the conditions outlined in Theorem~\ref{thm:main}. Recall from Section~\ref{sec:02} that there is no inclusion relationship between the set of trajectories outlined in Proposition~\ref{prop:Dubins3D} and $S(f_{3D}, t, \boldsymbol{\chi}_{3D, i})$. The only relationship is that they have H curves in common. Thus, application of Lemma~\ref{lemma:extendability} will be confined to H curves, and the arguments related to CSC and CCC curves must be addressed separately.

	From Proposition~\ref{prop:Dubins3D}, it follows that $bd(\mathcal{G}(f_{3D}, t, \boldsymbol{\chi}_{3D, i})) \subseteq A_1(t) \subseteq \mathcal{G}(f_{3D}, t, \boldsymbol{\chi}_{3D, i})$ and $bd(\mathcal{G}(-f_{3D}, T-t, \boldsymbol{\chi}_{3D, f})) \subseteq A_2(t) \subseteq \mathcal{G}(-f_{3D}, T-t, \boldsymbol{\chi}_{3D, f})$ for each $t$. Extendability of H curves follow from Lemma~\ref{lemma:extendability} since H curves are extremals. Extendability of CSC and CCC curves is trivial as discussed during the proof of the previous theorem. From extendability, we notice that $A_1(t)$ and $A_2(t)$ vary continuously in $t$ because the endpoints of all CSC, CCC, and H curves exist for all $t \in [0, T]$ and vary continuously. Proof of compactness of the set of endpoints of CSC and CCC curves follow the analogous steps in the previous theorem. Thus, it remains to study the limiting behavior of sequences of H curves. (i.e., to prove that the limit points of the endpoints of H curves belong to $A_1(t)$ and $A_2(t)$.)

	Consider a sequence of H curves denoted as $\{\boldsymbol{\chi}_k(t)\}_{k \in \mathbb{N}}$ which belongs to $S(f_{3D}, t, \boldsymbol{\chi}_{3D, i})$. From Lemma~\ref{lemma:extendability}, there exists a uniformly converging subsequence, denoted by the same notation $\boldsymbol{\chi}_k = (\boldsymbol{x}_k, \boldsymbol{y}_k)$, which its limit also belongs to $S(f_{3D}, t, \boldsymbol{\chi}_{3D, i})$. Let us denote such limit as $\overline{\boldsymbol{\chi}} = (\overline{\boldsymbol{x}}, \overline{\boldsymbol{y}})$. Since H curves are smooth, each $\boldsymbol{y}_k$ is differentiable everywhere. Moreover, curvature bound implies that each $\dot{\boldsymbol{y}}_k$ is uniformly bounded, and that all $\dot{\boldsymbol{y}}_k$ share the same Lipschitz constant of 1. This implies equicontinuity. Thus, $\{\boldsymbol{y}_k\}_{k \in \mathbb{N}}$ is a family of differentiable functions that uniformly converges to $\overline{\boldsymbol{y}}$, and have equicontinuous derivatives $\dot{\boldsymbol{y}}_k$. This implies that the uniform limit $\overline{\boldsymbol{y}}$ is differentiable and $\dot{\boldsymbol{y}}_k \rightarrow \dot{\overline{\boldsymbol{y}}}$ uniformly. Then $\dot{\overline{\boldsymbol{y}}}$ is continuous since uniform convergence preserves continuity. Now we recall that $S(f_{3D}, t, \boldsymbol{\chi}_{3D, i})$ consists certain concatenation of C and S components, and H curves. Thus, the only trajectories in $S(f_{3D}, t, \boldsymbol{\chi}_{3D, i})$ with continuous $\dot{\boldsymbol{y}}$ are H, C, and S. Since H curves have a constant curvature of 1, we further conclude that $\overline{\boldsymbol{\chi}}$ is either H or C. This proves that $\overline{\boldsymbol{\chi}}(t)$ is in $A_1(t)$, and compactness follows thereby. Thus, the choice of $A_1(t)$ and $A_2(t)$ satisfies the conditions outlined in Theorem~\ref{thm:main}. Then application of Theorem~\ref{thm:main} completes the proof. \qed
\end{pf*}
Similar to Theorem~\ref{thm:dubins2d}, Theorem~\ref{thm:dubins3d} facilitates trajectory synthesis by constructing, for each segment class pair, a residual map whose decision variables consist of the segment parameters (e.g., lengths and osculating plane angles for C and S segments, and the parameters of Eq.~\eqref{eq:torsionODE} for H segments) together with the concatenation time. The residual is defined as the mismatch between the two extremals at the concatenation point. The construction procedure is summarized in Algorithm~\ref{alg:3d_residual}. The existence of a solution to this formulation is equivalent to that of the original point-to-point steering problem. Nevertheless, it is important to note that the characterizations in Theorems~\ref{thm:dubins2d} and~\ref{thm:dubins3d} are based on the assumption that $\mathcal{B}(f, T, \boldsymbol{\chi}_i, \boldsymbol{\chi}_f)\neq\emptyset$. In general, there exist boundary and terminal time configurations for which no feasible solution exists and $\mathcal{B}(f, T, \boldsymbol{\chi}_i, \boldsymbol{\chi}_f)=\emptyset$. The formulations based on Theorems~\ref{thm:dubins2d} and~\ref{thm:dubins3d} do not provide a method to decide solution existence (i.e., whether $\mathcal{B}(f, T, \boldsymbol{\chi}_i, \boldsymbol{\chi}_f)=\emptyset$), as feasibility itself is a separate reachability question. Efficient feasibility tests for curvature-bounded paths in $\mathbb{R}^2$ are well studied in the literature~(\cite{patsko2003,patsko2022,Ding2019,chen2023,Rao2024}), and we leave the corresponding feasibility question in $\mathbb{R}^3$ as future work.

The congruence between the proofs of Theorem~\ref{thm:dubins2d} and Theorem~\ref{thm:dubins3d} is noticeable: the existing results on reachability are used to find the sets $A_1(t)$ and $A_2(t)$ in Theorem~\ref{thm:main}. This congruence arises because the arguments of Theorem~\ref{thm:main} do not rely on ingenious geometric arguments, but only on the fundamental connection between reachability and problem \textbf{P1}. While PMP was applied to prove Propositions~\ref{prop:Dubins2D} and~\ref{prop:Dubins3D} in~\cite{patsko2003} and~\cite{bae2025reachability}, the arguments of PMP do not explicitly appear during the proof.
\begin{algorithm}[htbp]
\scriptsize
\caption{Residual-Map Root-Finding for Curvature-Bounded Paths in $\mathbb{R}^3$}
\label{alg:3d_residual}
\begin{algorithmic}[1]
\Require Initial state $\boldsymbol{\chi}_{i}$, terminal state $\boldsymbol{\chi}_{f}$, terminal time $T$
\State Define segment classes $\mathcal{S} \gets \{\text{CSC, CCC, H}\}$
\State Initialize best residual norm $e^* \gets \infty$ and trajectory $\boldsymbol{\chi}^* \gets \emptyset$
\For{each combination $(C_1, C_2) \in \mathcal{S} \times \mathcal{S}$}
    \State \textbf{Define decision variables $\boldsymbol{\theta}$:}
    \State \quad For $C_k \in \{\text{CSC, CCC}\}$: $\boldsymbol{p}_k \gets [l_{k,1}, l_{k,2}, \boldsymbol{\phi}_k^T]^T$ \Comment{\newline \hspace*{5em} Lengths $l_{k,1}$, $l_{k,2}$ of the first two segments, and \newline \hspace*{5em} osculating plane angles $\boldsymbol{\phi}_k$ for each C segment (2 \newline \hspace*{5em} for CSC, 3 for CCC)}
    \State \quad For $C_k \in \{\text{H}\}$: $\boldsymbol{p}_k \gets [\varphi_k, \tau_{k}(0), \dot{\tau}_{k}(0), \zeta_k]^T$ \Comment{Initial \newline \hspace*{5em} osculating plane angle $\varphi_k$ and params for Eq.~\eqref{eq:torsionODE}}
    \State \quad $\boldsymbol{\theta} \gets [\boldsymbol{p}_1^T, \boldsymbol{p}_2^T, t_c]^T$ \Comment{Concatenation time $t_c \in [0, T]$}
    
    \Function{Residual}{$\boldsymbol{\theta}$}
        
        \State \textbf{Forward integration for $C_1$ over $[0, t_c]$:}
        \If{$C_1 == \text{H}$}
            \State Obtain $\boldsymbol{\chi}_1(t_c)$ by integrating the Frenet-Serret \newline \hspace*{5em} equations from $\boldsymbol{\chi}_{i}$ with $\kappa \equiv 1$ and torsion $\tau(t)$ \newline \hspace*{5em} parametrized by Eq.~\eqref{eq:torsionODE} and $\boldsymbol{p}_1$
        \Else
            \State Construct $\boldsymbol{\chi}_1(t_c)$ geometrically from $\boldsymbol{\chi}_{i}$ using $\boldsymbol{p}_1$
        \EndIf
        
        \State \textbf{Backward integration for $C_2$ over $[t_c, T]$:}
        \If{$C_2 == \text{H}$}
            \State Obtain $\boldsymbol{\chi}_2(t_c)$ by integrating the Frenet-Serret \newline \hspace*{5em} equations backward from $\boldsymbol{\chi}_{f}$ with $\kappa \equiv 1$ and \newline \hspace*{5em} torsion $\tau(t)$ parametrized by Eq.~\eqref{eq:torsionODE} and $\boldsymbol{p}_2$
        \Else
            \State Construct $\boldsymbol{\chi}_2(t_c)$ geometrically from $\boldsymbol{\chi}_{f}$ using $\boldsymbol{p}_2$
        \EndIf
        
        \State \Return $\boldsymbol{\chi}_1(t_c) - \boldsymbol{\chi}_2(t_c)$
    \EndFunction
    
    \State Solve \textsc{Residual}($\boldsymbol{\theta}$) = $\boldsymbol{0}$ using a nonlinear equation solver
    \If{Solution converges with residual norm $e < e^*$}
        \State $e^* \gets e$
        \State $\boldsymbol{\chi}^* \gets [\boldsymbol{\chi}_1([0, t_c]), \boldsymbol{\chi}_2([t_c, T])]$
    \EndIf
\EndFor
\State \Return $\boldsymbol{\chi}^*$
\end{algorithmic}
\end{algorithm}
\begin{rem}
	The structure of extremals of curvature-bounded paths in general dimensions (and even on spheres and hyperbolic spaces) is obtained by embedding the three-dimensional extremals into a higher-dimensional ambient space in~\cite{Mittenhuber1998}. Hence, Theorem~\ref{thm:main} is expected to extend to general dimensions in the same spirit as Theorem~\ref{thm:dubins3d}. A rigorous treatment, however, requires a reformulation of the dynamics on the Lie group $SE(n)$ of Euclidean motions of $\mathbb{R}^n$, which we leave for future work.
\end{rem}
%
\section{Numerical Demonstrations}\label{sec:05}
We present numerical demonstrations to validate the efficacy of the TPBVP formulation (problem \textbf{P2}) for point-to-point steering problem within the setting of Lemma~\ref{lemma:extendability} and where Pontryagin extremals are computable. We first consider a controlled Van der Pol system as a nonlinear example to illustrate applicability beyond curvature-bounded path dynamics. We then proceed to curvature-bounded paths in $\mathbb{R}^2$ and $\mathbb{R}^3$ with Theorems~\ref{thm:dubins2d} and~\ref{thm:dubins3d} as practical application classes where extremals admit compact descriptions. For Dubins car (curvature-bounded paths in $\mathbb{R}^2$), the point-to-point steering problem admits well-established geometric constructions in the literature~(\cite{Ding2019,chen2023,Rao2024}). We therefore include a direct comparison with~\cite{chen2023} as a reference point, to confirm that the algebraic residual-map formulation implied by Theorem~\ref{thm:dubins2d} reproduces established solutions. Importantly, we note that the main emphasis is on demonstrating that the same theoretical framework applies beyond geometry-specific constructions, as illustrated by the nonlinear Van der Pol example and the curvature-bounded paths in $\mathbb{R}^3$. For the latter, we demonstrate that analogous trajectory synthesis results can be obtained through Theorem~\ref{thm:dubins3d} under the reachability assumption, in a setting where directly comparable geometric benchmarks are not available.

Generally speaking, TPBVPs (problem \textbf{P2}) and the residual equations associated with Theorems~\ref{thm:dubins2d} and~\ref{thm:dubins3d} may admit multiple solutions. In the latter, additional multiplicity arises because different segment combinations induce different residual maps. Therefore, multiple solutions may be obtained even when the same initial guess is used within each segment combination. When multiple solutions are encountered in the experiments, we report the one with the smallest solver residual. All code used is available online\footnote{\href{https://github.com/johnbae1901/Fixed-Time-State-Transfer}{https://github.com/johnbae1901/Fixed-Time-State-Transfer}}.
\subsection{Controlled Van der Pol system}
We apply the TPBVP formulation (problem \textbf{P2}) to a controlled Van der Pol system below with a control restraint set $\Omega = [-1, 1]$.
\begin{equation} \label{eq:vdp_dynamics}
	\begin{split}
		& \dot{x}_1 = x_2 \\
		& \dot{x}_2 = \left( 1 - x_1^2 \right) x_2 - x_1 + u
	\end{split}
\end{equation}
By denoting costate variable as $\boldsymbol{p} = (p_1, p_2)$, it is easy to verify that the maximizer of Hamiltonian of this system is given as $u^* = sgn(p_2)$ and that no singular arc exists. Since problem P2 is a TPBVP with an additional unknown parameter $\tau$, we adopt a collocation-based BVP solver \emph{BoundaryValueDiffEq.jl}~(\cite{rackauckas2017}) in SciML, which is typically more robust than shooting method for addressing nonlinear dynamics such as Van der Pol. As in other dynamics, feasibility of a given terminal time and boundary values is a reachability question and is outside the scope of this paper; here we select representative feasible instances for demonstration.

We present the results for two sets of boundary configurations in Table~\ref{tab:vdp_summary}, which summarizes both the prescribed problem instances and the corresponding computed parameters. For each configuration, we consider two variants: one with fully constrained terminal state, and the other with partially constrained terminal state where $x_2(T)$ is free. When $x_2(T)$ is free, it is treated as an additional unknown parameter of the TPBVP and obtained as part of the solution (cf. Remark~\ref{rem:unconstrained}). As illustrated in Figure~\ref{fig:VdP_results}, the solution trajectory is successfully constructed by joining two extremals originating from $\boldsymbol{\chi}_i$ and $\boldsymbol{\chi}_f$ respectively. The corresponding control inputs exhibit a bang-bang switching structure consistent with the maximization condition $u^* = sgn(p_2)$.
\begin{table}[t]
    \centering
    \resizebox{\columnwidth}{!}{%
        \tiny
        \setlength{\tabcolsep}{1pt}
        \begin{tabular}{@{} c | c c c | c c c | c c c | c c c @{}}
            \toprule
            \multirow{3}{*}{\raisebox{-3.2ex}{\textbf{Item}}} &
            \multicolumn{6}{c|}{\textbf{Inst. \#1} ($\boldsymbol{\chi}_f=(0.0,\cdot)$)} &
            \multicolumn{6}{c}{\textbf{Inst. \#2} ($\boldsymbol{\chi}_f=(0.6,\cdot)$)} \\
            \cmidrule(lr){2-7} \cmidrule(lr){8-13}
            &
            \multicolumn{3}{c|}{\textbf{(a) fixed $x_2(T)$}} &
            \multicolumn{3}{c|}{\textbf{(b) free $x_2(T)$}} &
            \multicolumn{3}{c|}{\textbf{(a) fixed $x_2(T)$}} &
            \multicolumn{3}{c}{\textbf{(b) free $x_2(T)$}} \\
            \cmidrule(lr){2-4} \cmidrule(lr){5-7} \cmidrule(lr){8-10} \cmidrule(lr){11-13}
            & $x_{i,1}$ & $x_{i,2}$ & $T$ & $x_{i,1}$ & $x_{i,2}$ & $T$ & $x_{i,1}$ & $x_{i,2}$ & $T$ & $x_{i,1}$ & $x_{i,2}$ & $T$ \\
            \midrule
            $\boldsymbol{\chi}_i$ &
            $2.0$ & $2.0$ & $4.0$ &
            $2.0$ & $2.0$ & $4.0$ &
            $2.0$ & $2.0$ & $4.0$ &
            $2.0$ & $2.0$ & $4.0$ \\

            $\boldsymbol{\chi}_f$ &
            $0.0$ & $0.0$ & -- &
            $0.0$ & \textit{free} & -- &
            $0.6$ & $-0.9$ & -- &
            $0.6$ & \textit{free} & -- \\

            $x_2(T)$ &
            $0.0$ & -- & -- &
            $-0.3585$ & -- & -- &
            $-0.9$ & -- & -- &
            $-0.3701$ & -- & -- \\

            $\tau$ &
            $0.7343$ & -- & -- &
            $0.6313$ & -- & -- &
            $0.4453$ & -- & -- &
            $0.4077$ & -- & -- \\

            $\tau T$ &
            $2.9373$ & -- & -- &
            $2.5254$ & -- & -- &
            $1.7812$ & -- & -- &
            $1.6307$ & -- & -- \\
            \bottomrule
        \end{tabular}
    }
    \caption{Summary of the configurations and results of Van der Pol system simulation. Two sets of boundary conditions are considered, each with two terminal scenarios: (a) fixed and (b) free $x_2(T)$. The prescribed boundary values ($\boldsymbol{\chi}_i, \boldsymbol{\chi}_f$) as well as the computed parameters ($\tau$, concatenation time $\tau T$, and the obtained $x_2(T)$ for the free cases) are listed.}
    \label{tab:vdp_summary}
\end{table}
\newlength{\imgwVDP}
\setlength{\imgwVDP}{\dimexpr(\columnwidth-2\tabcolsep)/2\relax}
\begin{figure}[htbp]
  \centering
  \setlength{\tabcolsep}{2pt}

  \begin{tabular}{@{}cc@{}}

    \multicolumn{2}{@{}l@{}}{\scriptsize\textbf{Inst. \#1 (fixed $x_2(T)$)}}\\[-0.5mm]
    \includegraphics[width=\imgwVDP, trim=50 0 0 0]{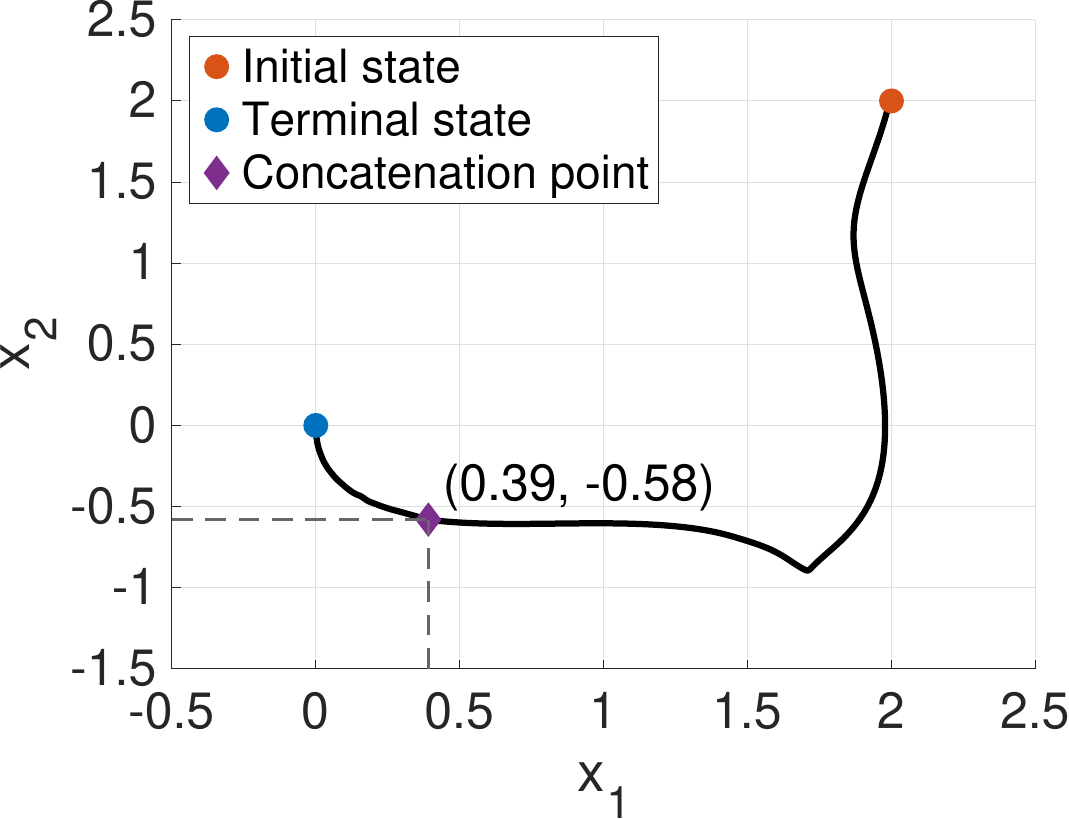} &
    \includegraphics[width=\imgwVDP, trim=50 0 0 0]{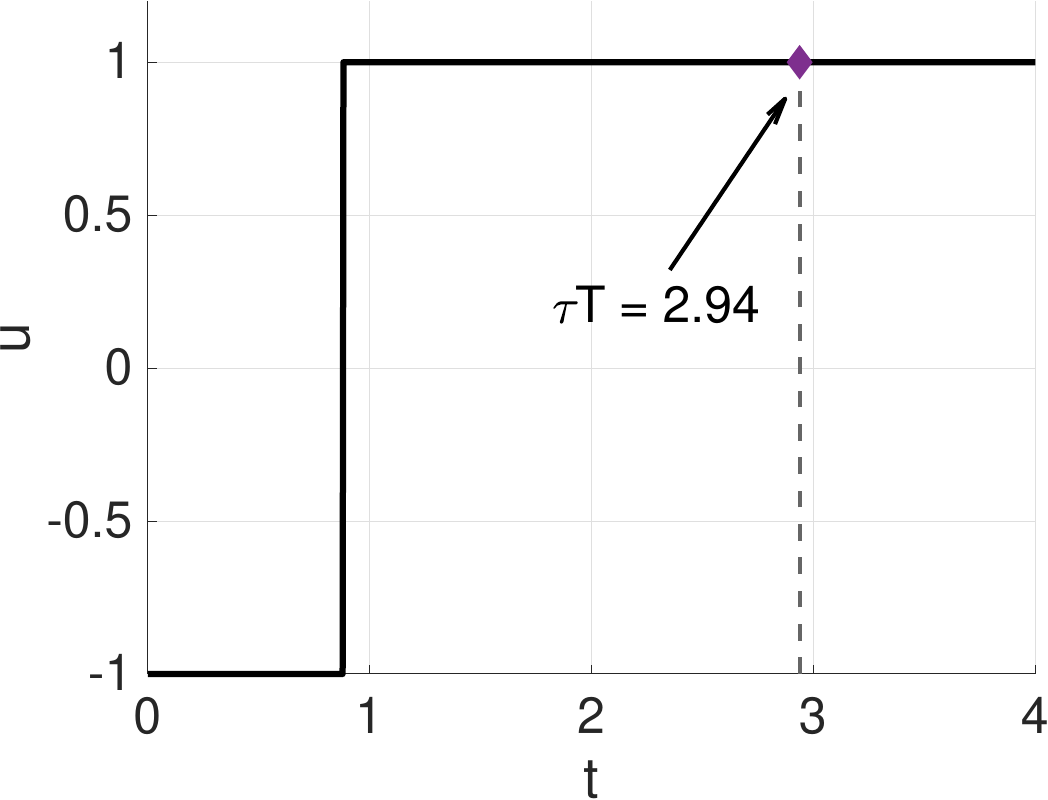} \\[-0.7mm]

    \multicolumn{2}{@{}l@{}}{\scriptsize\textbf{Inst. \#1 (free $x_2(T)$)}}\\[-0.5mm]
    \includegraphics[width=\imgwVDP, trim=50 0 0 0]{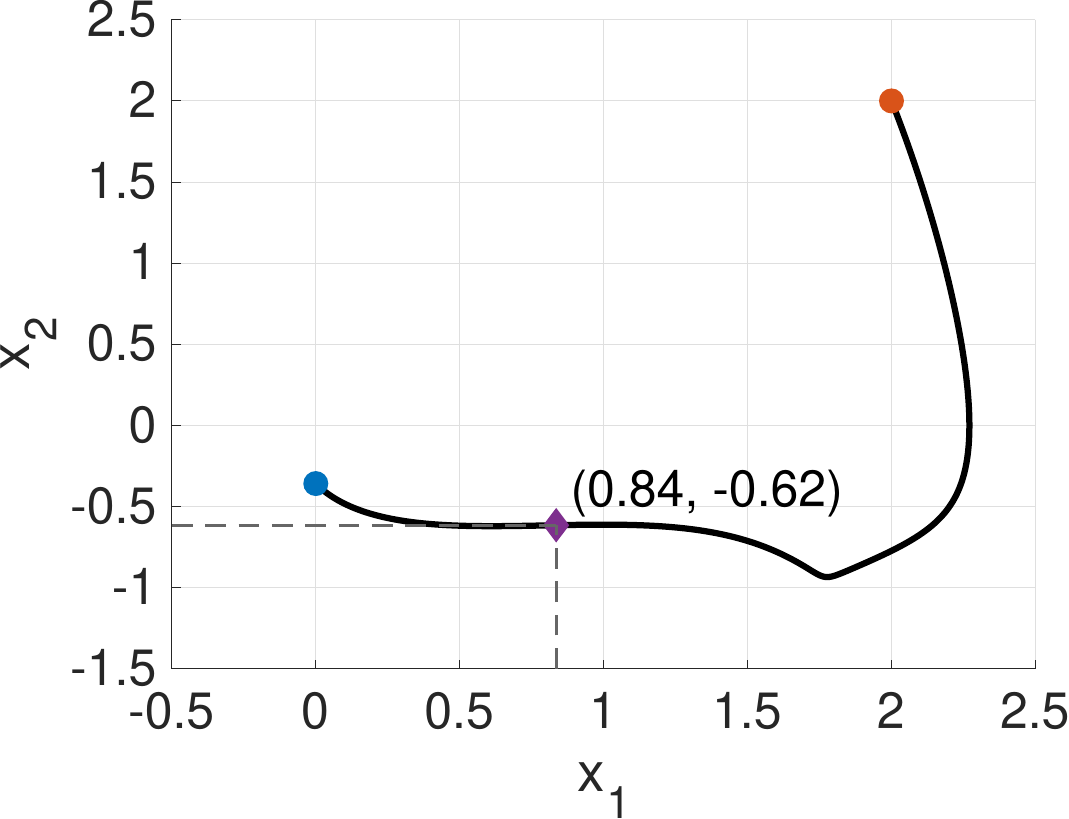} &
    \includegraphics[width=\imgwVDP, trim=50 0 0 0]{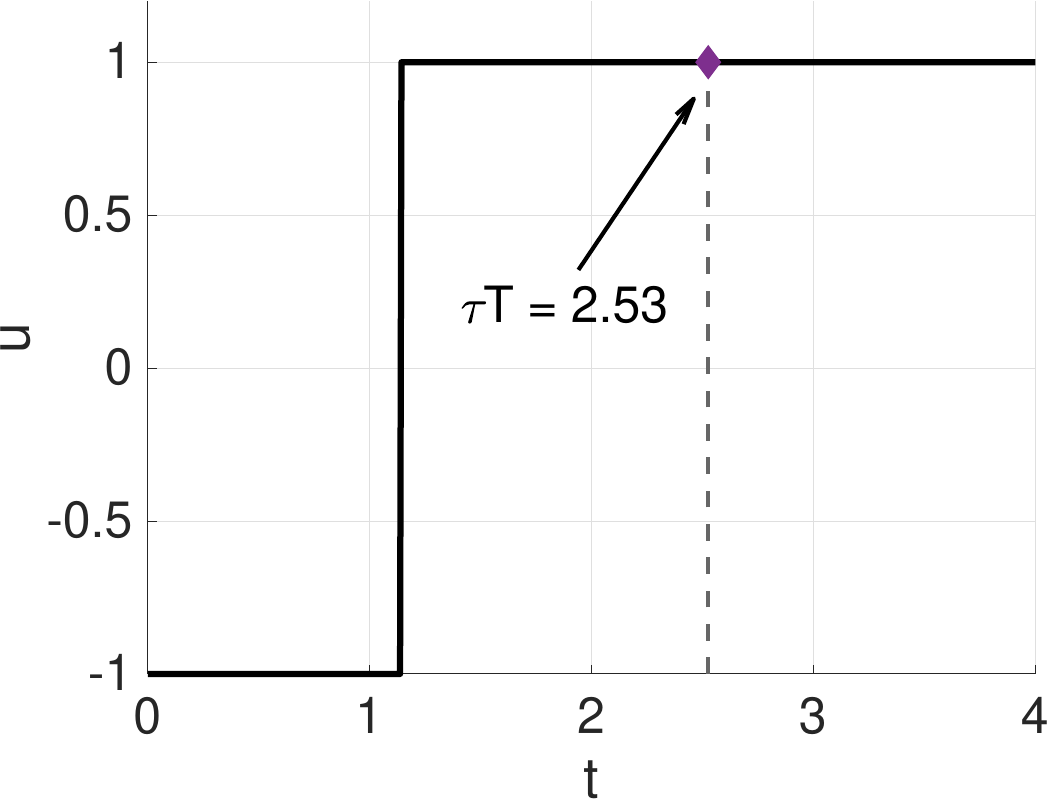} \\[-0.7mm]

    \multicolumn{2}{@{}l@{}}{\scriptsize\textbf{Inst. \#2 (fixed $x_2(T)$)}}\\[-0.5mm]
    \includegraphics[width=\imgwVDP, trim=50 0 0 0]{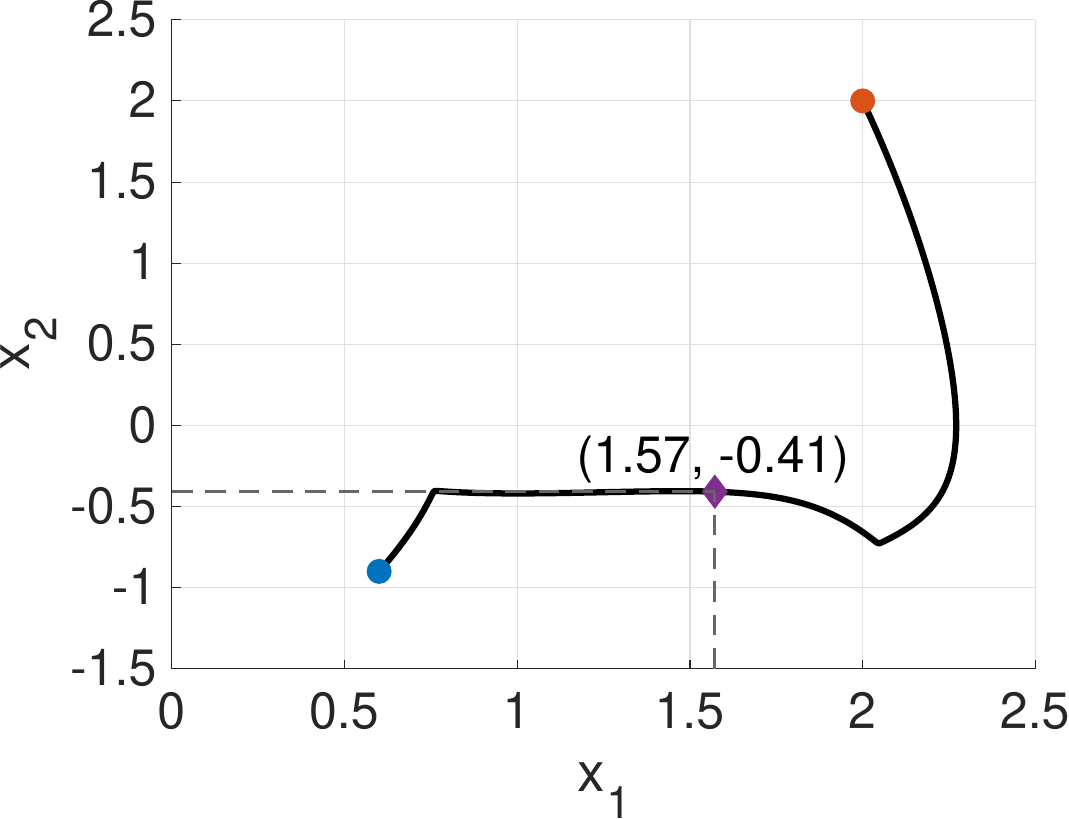} &
    \includegraphics[width=\imgwVDP, trim=50 0 0 0]{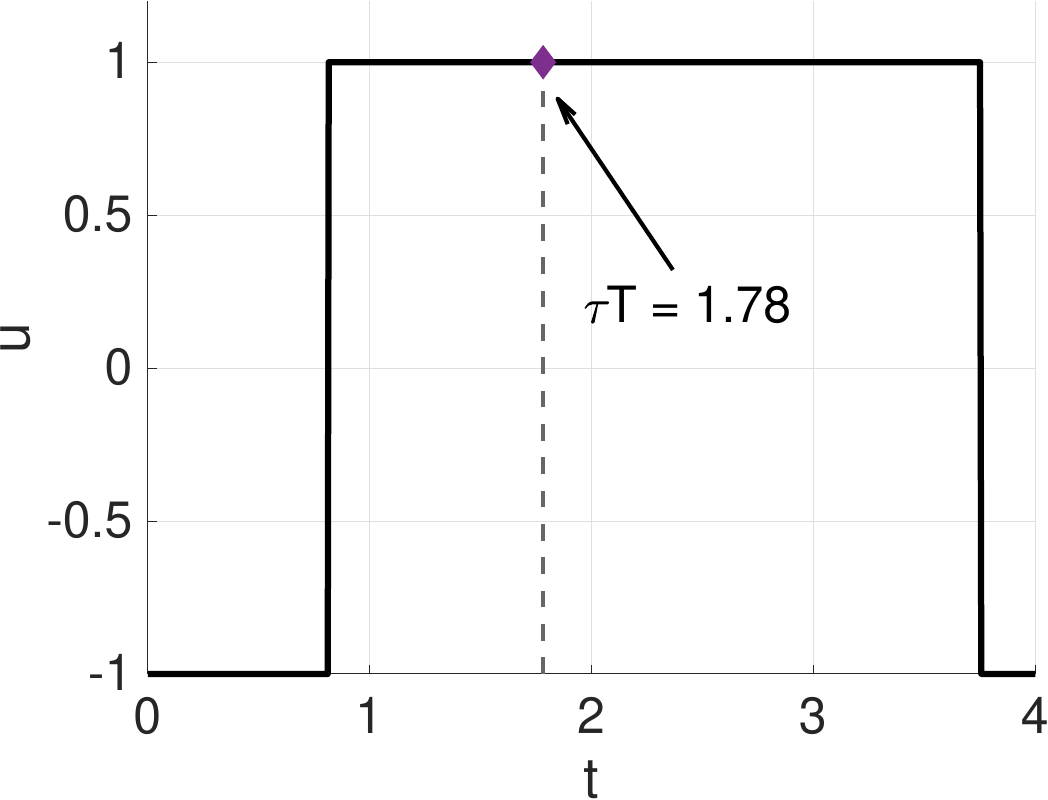} \\[-0.7mm]

    \multicolumn{2}{@{}l@{}}{\scriptsize\textbf{Inst. \#2 (free $x_2(T)$)}}\\[-0.5mm]
    \includegraphics[width=\imgwVDP, trim=50 0 0 0]{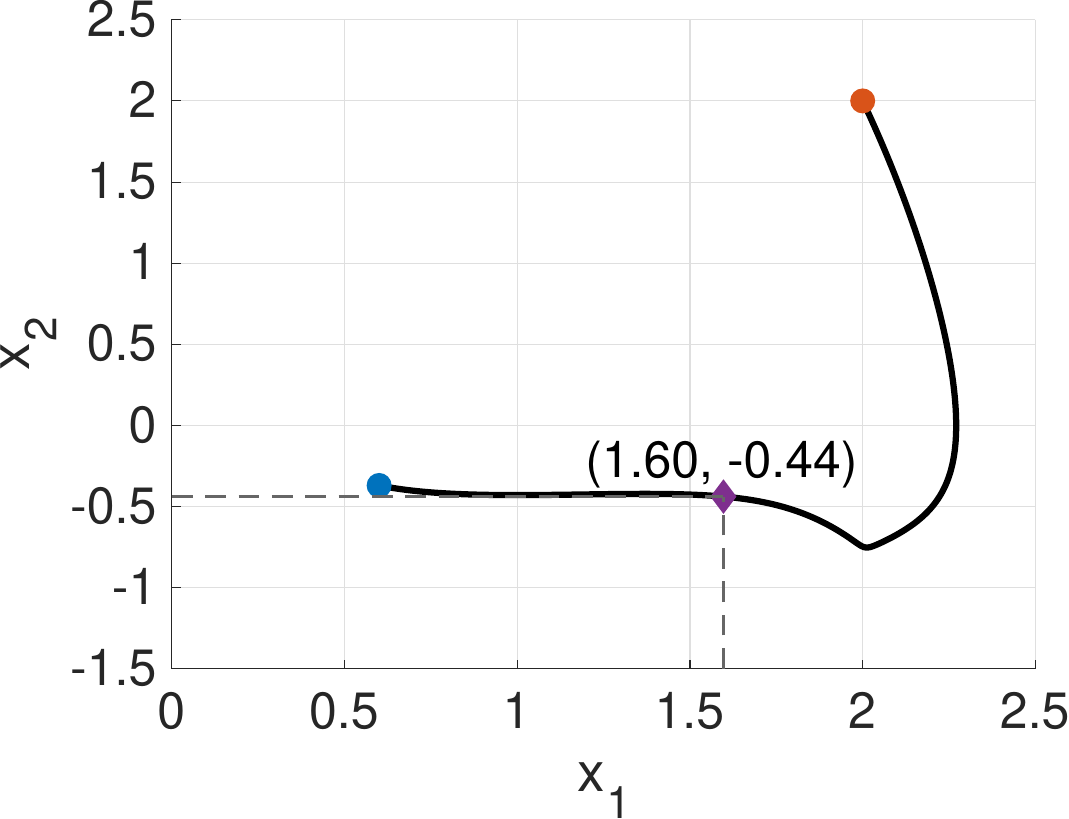} &
    \includegraphics[width=\imgwVDP, trim=50 0 0 0]{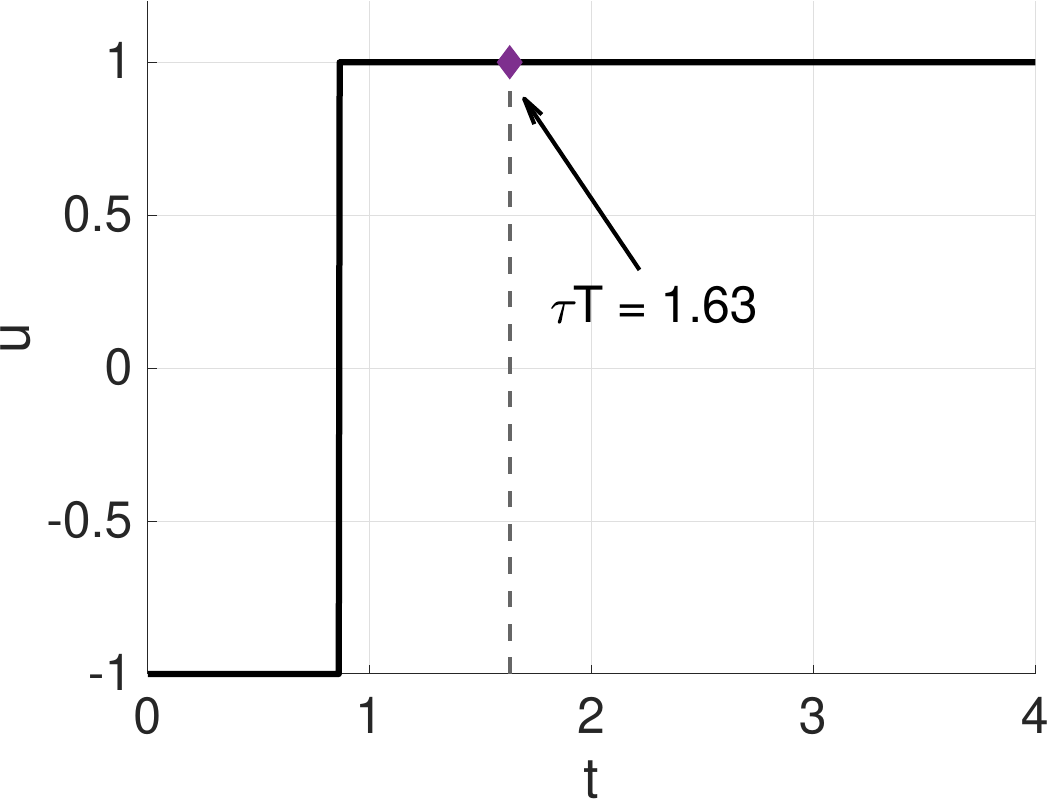} \\[-0.7mm]

    \scriptsize \textbf{Trajectory $(x_1(t), x_2(t))$} & \scriptsize \textbf{Control input $u(t)$} \\[-0.1mm]
  \end{tabular}
  \caption{Numerical solutions to the point-to-point steering problem of controlled Van der Pol system obtained by solving Problem \textbf{P2}. The rows correspond to the four scenarios listed in Table~\ref{tab:vdp_summary}. \textbf{Left:} State-space trajectories in the $(x_1, x_2)$ plane. Orange and blue markers denote the initial ($\boldsymbol{\chi}_i$) and terminal ($\boldsymbol{\chi}_f$) states, respectively. The purple diamond indicates the concatenation point where the two extremals join. \textbf{Right:} Time history of the control input $u(t)$. The vertical dashed line marks the concatenation time $t = \tau T$.}
  \label{fig:VdP_results}
\end{figure}
\subsection{Curvature-bounded paths in $\mathbb{R}^2$} \label{subsec:Dubins2D}
Theorem~\ref{thm:dubins2d} yields an algebraic residual map of the boundary conditions whose decision variables are the lengths of the C and S components (we omit the explicit expression for brevity). Unlike geometric elongation approaches, this enables the use of off-the-shelf nonlinear equation solvers to synthesize a curvature-bounded trajectory of prescribed length and boundary conditions. For demonstration, we consider a trajectory planning problem in which multiple Dubins cars are required to form a prescribed triangular formation simultaneously. The flexibility of the proposed formulation is that two scenarios of fixed and free terminal heading angles can be addressed by simply including it as an additional decision variable in the residual map (see Remark~\ref{rem:unconstrained}). We note that the free terminal heading case has also been addressed geometrically in~\cite{Ding2019}. We assume that terminal time is given appropriately so that the problem is feasible; testing feasibility of a given terminal time, or finding a feasible time given boundary conditions is, in essence, a reachability analysis problem which is addressed in~\cite{patsko2003,patsko2022,chen2023,Rao2024} for Dubins car. Specifically, the results in~\cite{chen2023} provide the explicit characterization of connected components of the set of feasible lengths for given boundary values. Once a feasible terminal time is given, the residual map constructed in Theorem~\ref{thm:dubins2d} admits at least one zero, and each zero corresponds to an admissible trajectory of the prescribed length.

\begin{figure}[htbp]
	\centering
	\includegraphics[scale=0.3, trim=40 0 0 0]{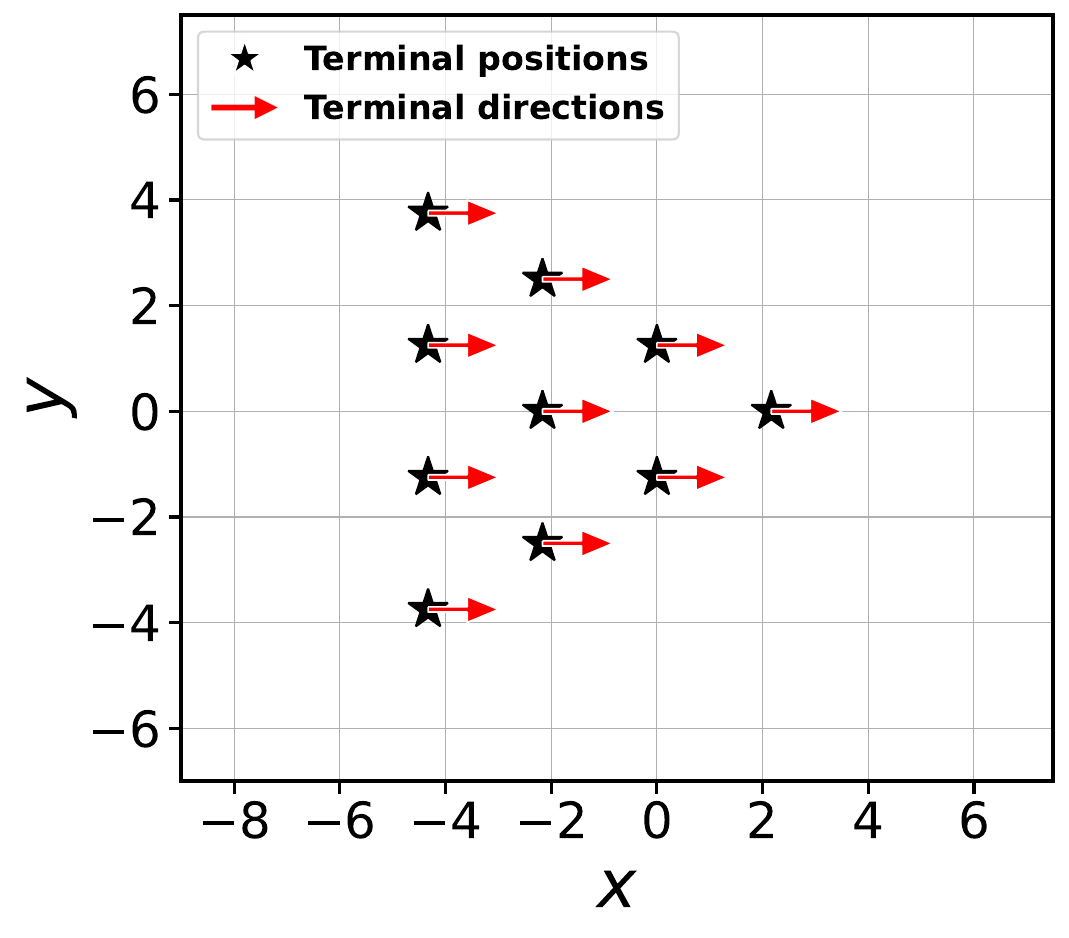}
	\caption{Setup of formation planning problem of Dubins car.}
	\label{fig:Dubins2D_formation_setup}
\end{figure}

Numerical experiment is conducted with 10 Dubins cars, required to form an equilateral triangle with distances of $2.5$ between each other, as illustrated in Figure~\ref{fig:Dubins2D_formation_setup}. Initial states are drawn uniformly randomly from the set $[-5, 5]^2 \times [0, 2\pi]$. We visualize the results for three representative instances generated from these random initial states. For each instance, a common terminal time $T$ for all 10 vehicles is set to the minimum feasible value extracted from the aforementioned connected components characterized in~\cite{chen2023}. The exact boundary values and the corresponding terminal times are available in the provided code repository. For each case, two scenarios are considered\textemdash{}one with fixed terminal heading angles, and the other with free terminal heading angles for each vehicle, using the same terminal time $T$. We use \textsc{Matlab} \textit{fsolve} to find a zero of the residual function. Initial guesses are set by assigning equal length of $T/6$ to each C and S component.

The obtained results for the two scenarios are illustrated in Figures~\ref{fig:Dubins2D_formation_fixed} and~\ref{fig:Dubins2D_formation_free} respectively. Figure~\ref{fig:Dubins2D_formation_fixed} (fixed terminal heading angles) shows that all vehicles reach the desired triangular formation at the prescribed time using the structures characterized in Theorem~\ref{thm:dubins2d}; concatenation of two trajectories of CSC, CCC, or their subsegments. Figure~\ref{fig:Dubins2D_formation_free} (free terminal heading angles) further demonstrates that the same framework accommodates free terminal headings by treating $\gamma(T)$ as an additional decision variable, as described in Remark~\ref{rem:unconstrained}. The resulting terminal heading angles for the free scenario, visualized as red arrows in Figure~\ref{fig:Dubins2D_formation_free}, are detailed in the provided code repository.
\newlength{\imgw}
\setlength{\imgw}{\dimexpr(\columnwidth-2\tabcolsep)/2\relax}
\begin{figure}[htbp]
  \centering
  \setlength{\tabcolsep}{2pt}
  \begin{tabular}{@{}cc@{}}
    \scriptsize \textbf{\cite{chen2023}} & \scriptsize \textbf{Theorem~\ref{thm:dubins2d}} \\[-0.1mm]

    \multicolumn{2}{@{}l@{}}{\scriptsize\textbf{Inst. \#1}}\\[-0.5mm]
    \includegraphics[width=\imgw, trim=20 0 0 0]{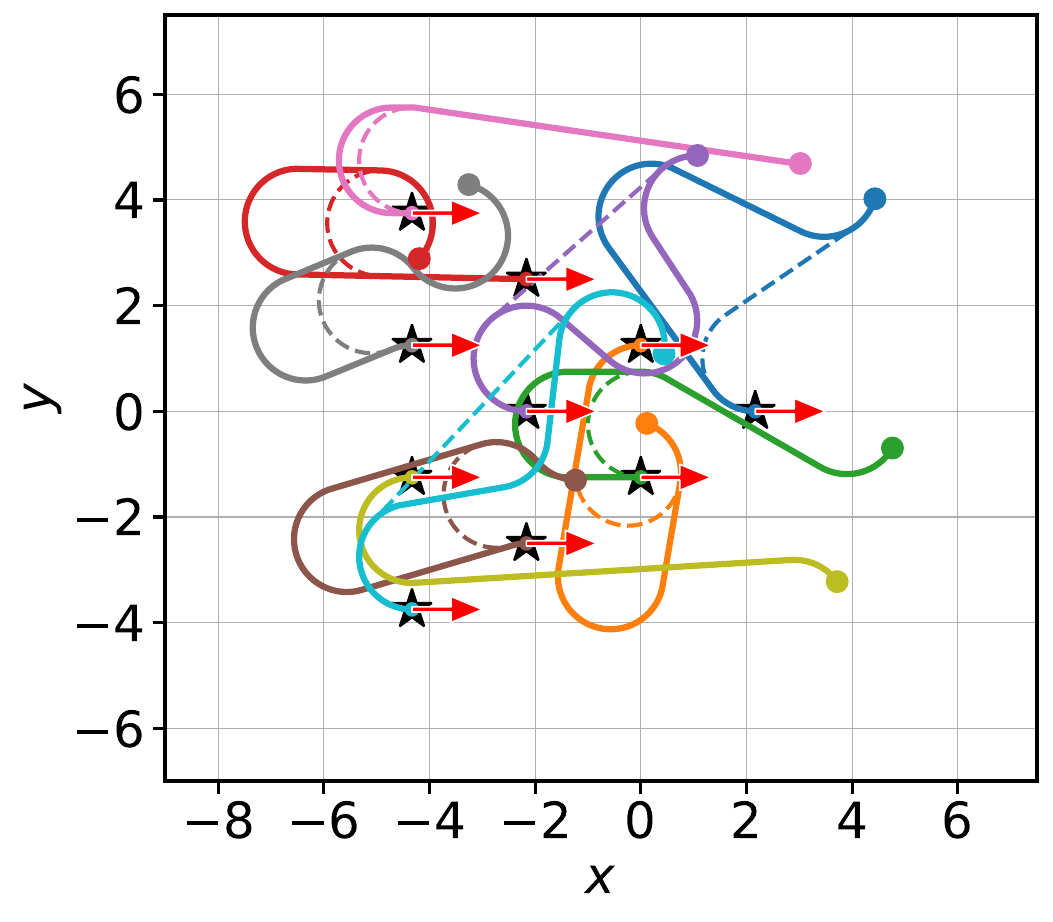} &
    \includegraphics[width=\imgw, trim=20 0 0 0]{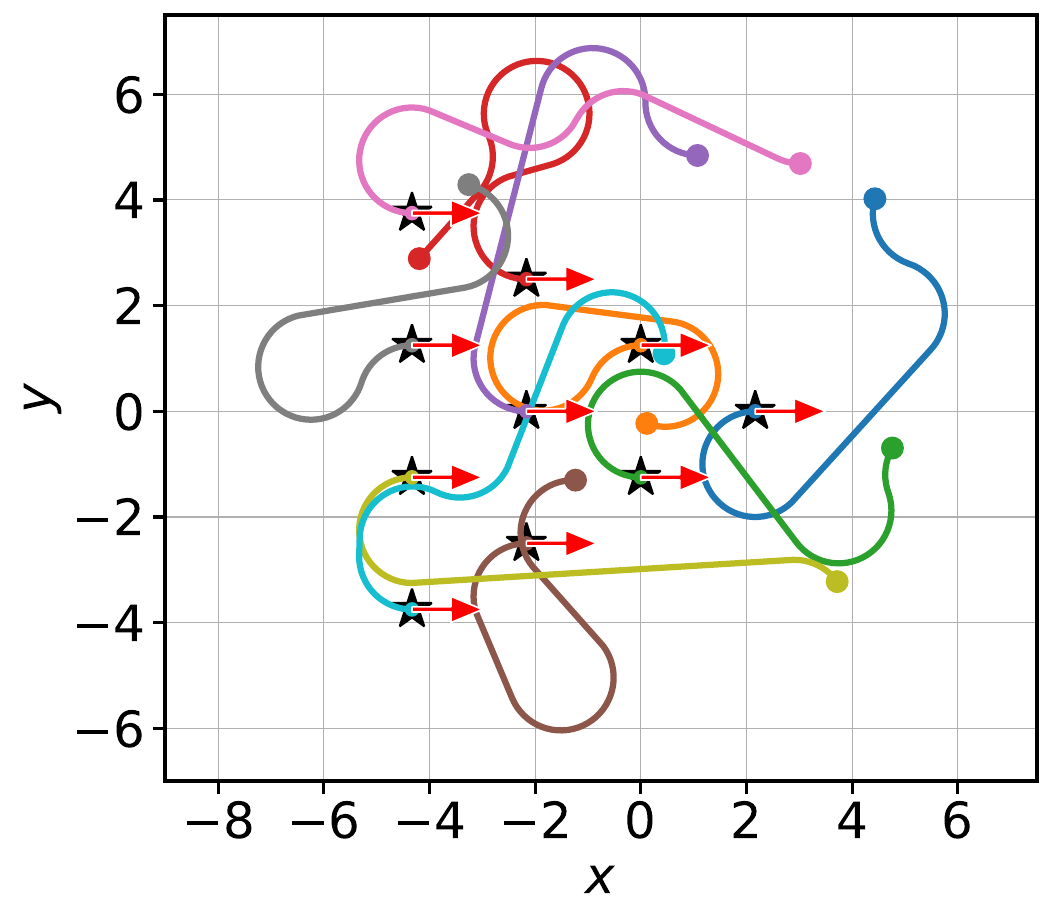} \\[-0.7mm]

    \multicolumn{2}{@{}l@{}}{\scriptsize\textbf{Inst. \#2}}\\[-0.5mm]
    \includegraphics[width=\imgw, trim=20 0 0 0]{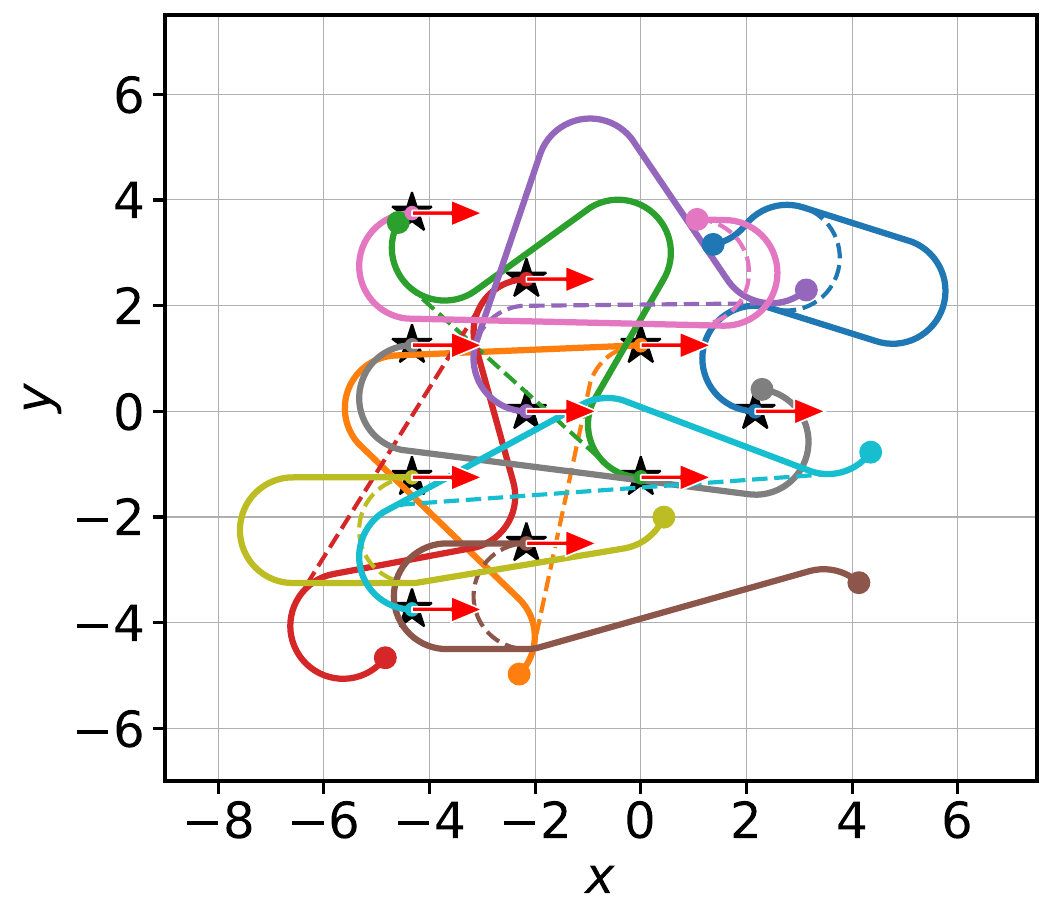} &
    \includegraphics[width=\imgw, trim=20 0 0 0]{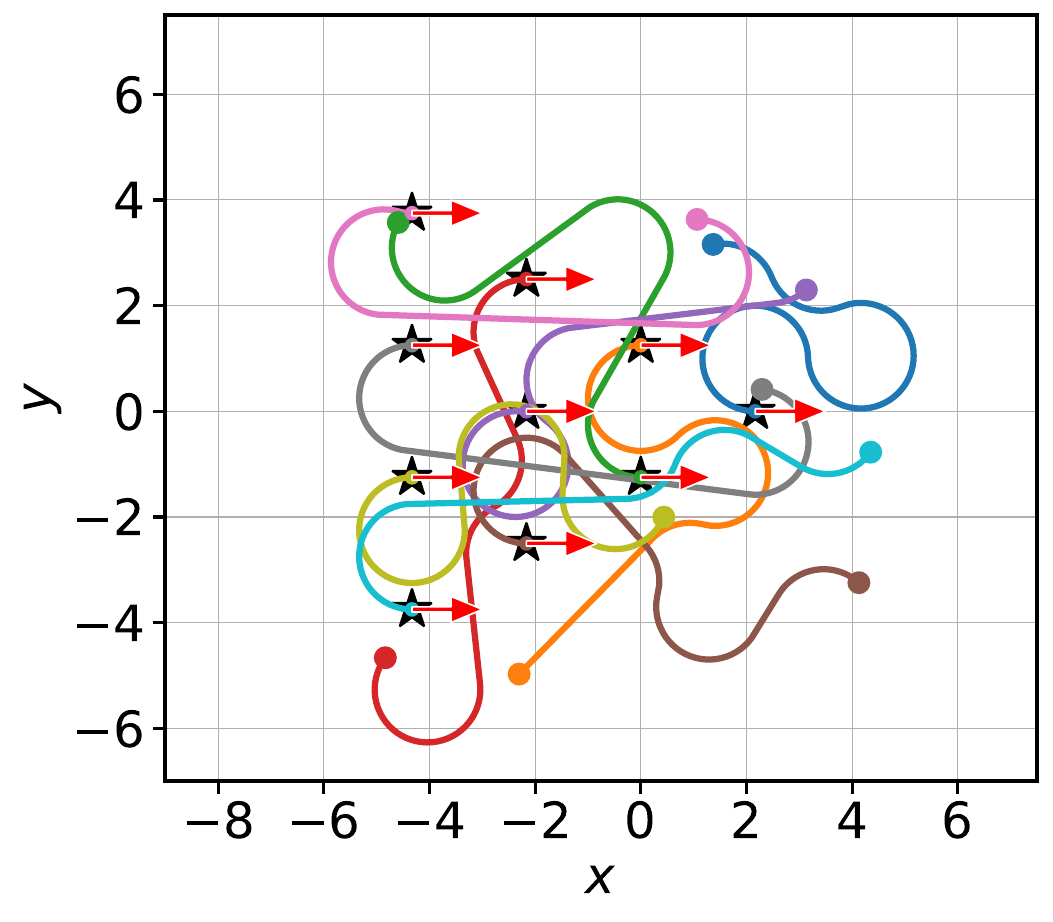} \\[-0.7mm]

    \multicolumn{2}{@{}l@{}}{\scriptsize\textbf{Inst. \#3}}\\[-0.5mm]
    \includegraphics[width=\imgw, trim=20 0 0 0]{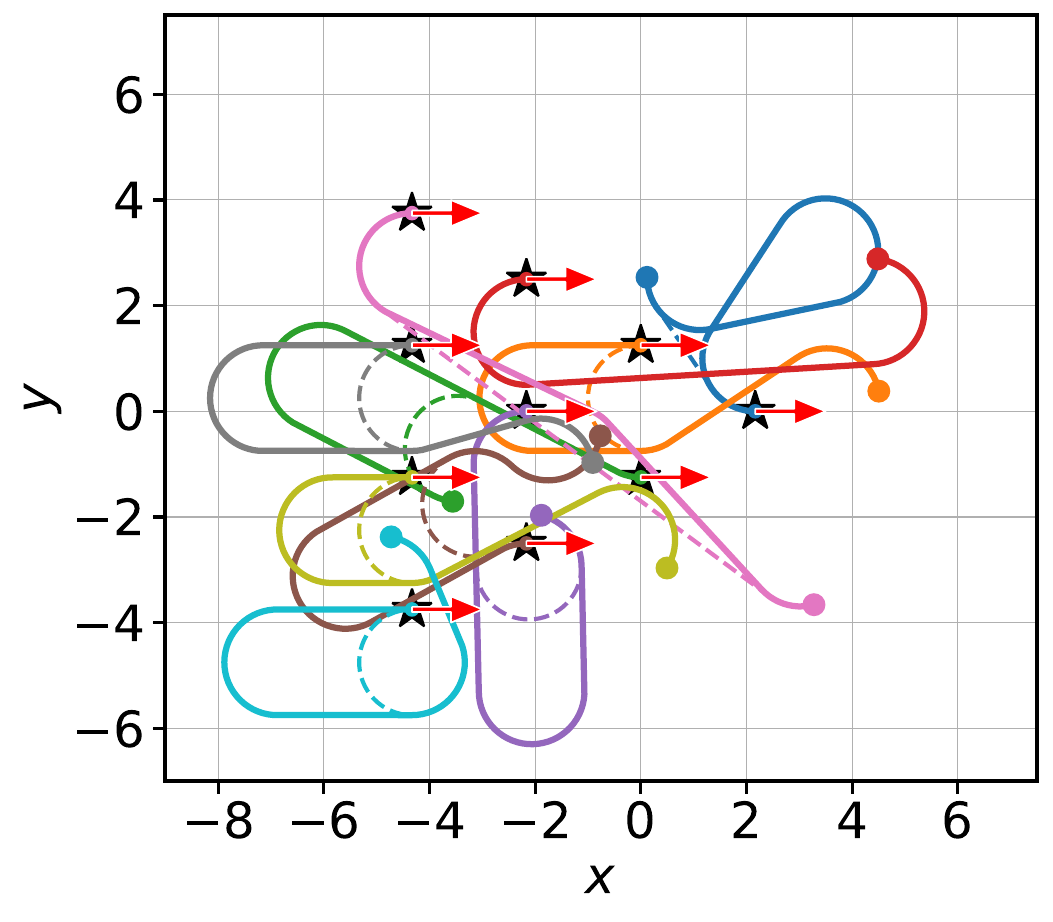} &
    \includegraphics[width=\imgw, trim=20 0 0 0]{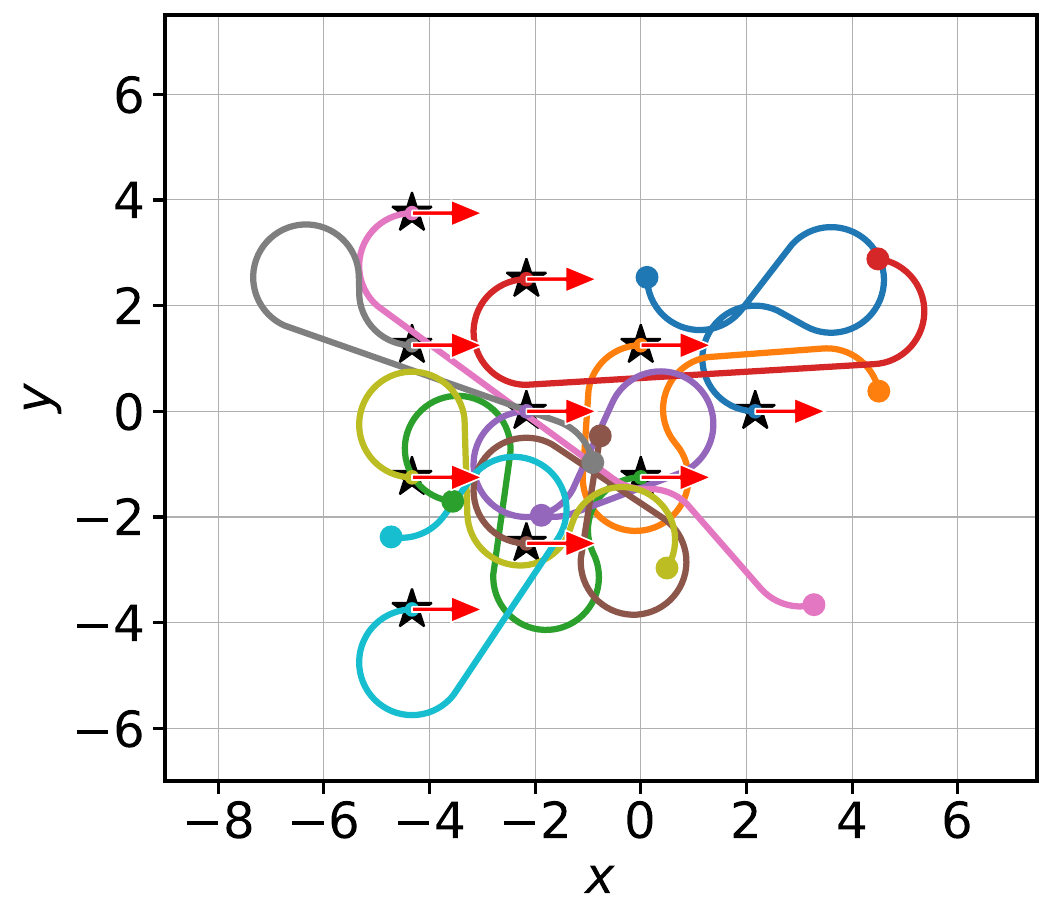} \\
  \end{tabular}
  \caption{Formation planning with fixed terminal heading angles for the three representative instances generated from random initial states. \textbf{Left}: geometric construction in~\cite{chen2023}. Dotted lines indicate the minimum-length paths used as references in the geometric construction of~\cite{chen2023}. \textbf{Right}: trajectories obtained from the residual map implied by Theorem~\ref{thm:dubins2d}.}
  \label{fig:Dubins2D_formation_fixed}
\end{figure}
\newlength{\imgwFreeTwo}
\setlength{\imgwFreeTwo}{\dimexpr(\columnwidth-2\tabcolsep)/2\relax}
\begin{figure}[htbp]
  \centering
  \setlength{\tabcolsep}{2pt}
  \renewcommand{\arraystretch}{1.0}

  \begin{tabular}{@{}cc@{}}
    \makebox[\imgwFreeTwo][c]{\scriptsize\textbf{Inst. \#1}} &
    \makebox[\imgwFreeTwo][c]{\scriptsize\textbf{Inst. \#2}} \\[-0.5mm]
    \includegraphics[width=0.9\imgwFreeTwo, trim=60 0 0 0]{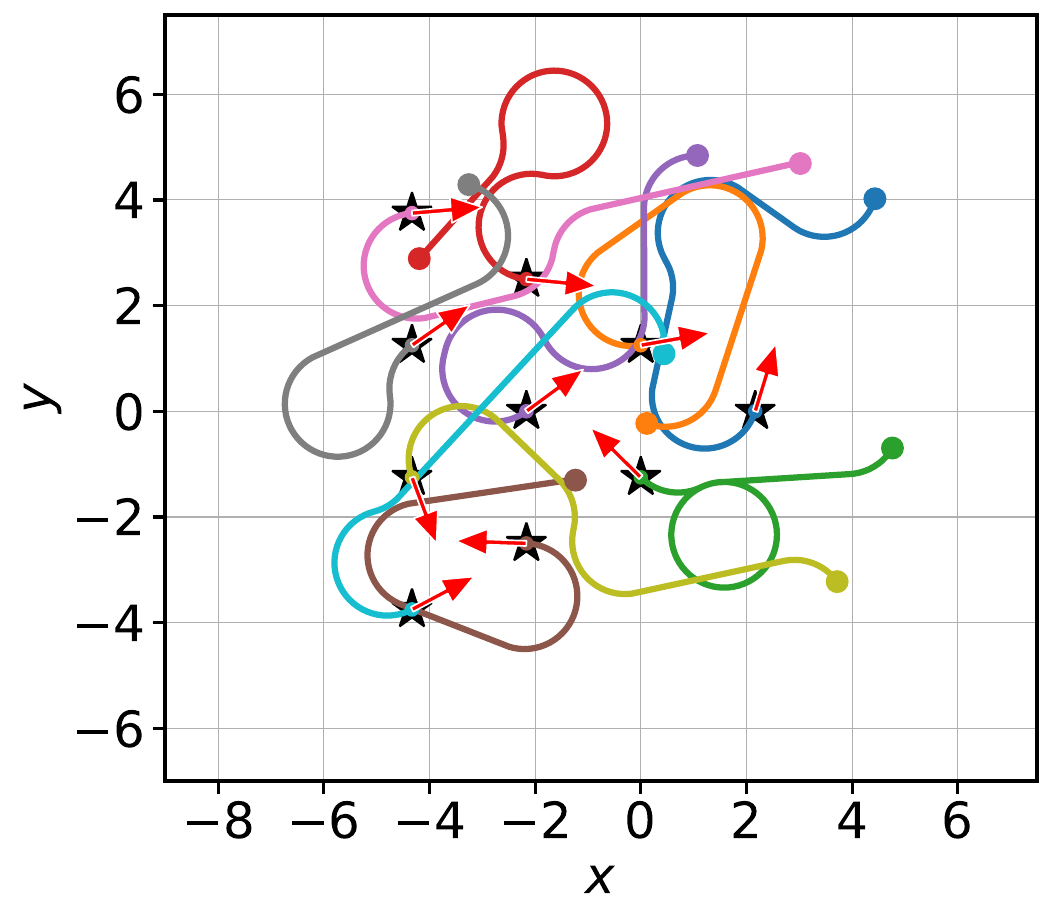} &
    \includegraphics[width=0.9\imgwFreeTwo, trim=60 0 0 0]{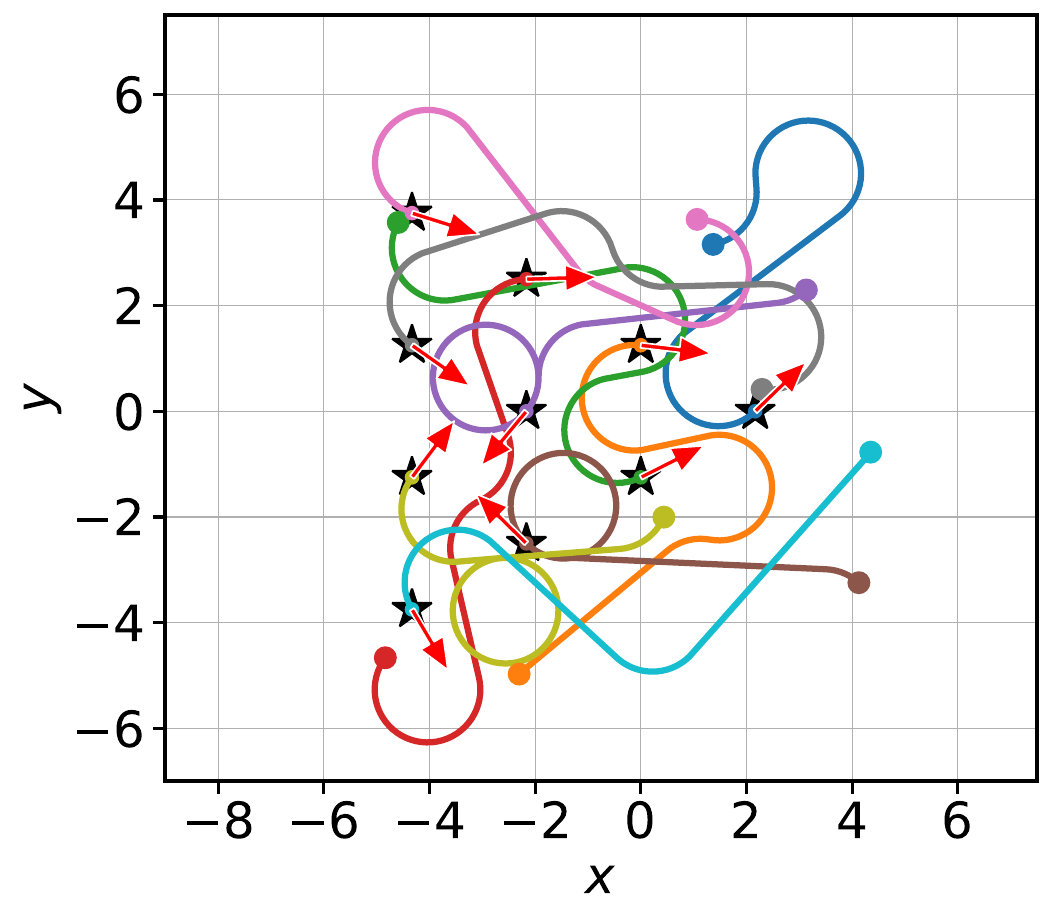} \\[-0.6mm]

    \multicolumn{2}{c}{\makebox[0.9\imgwFreeTwo][c]{\scriptsize\textbf{Inst. \#3}}}\\[-0.5mm]
    \multicolumn{2}{c}{\includegraphics[width=0.98\imgwFreeTwo, trim=60 0 0 0]{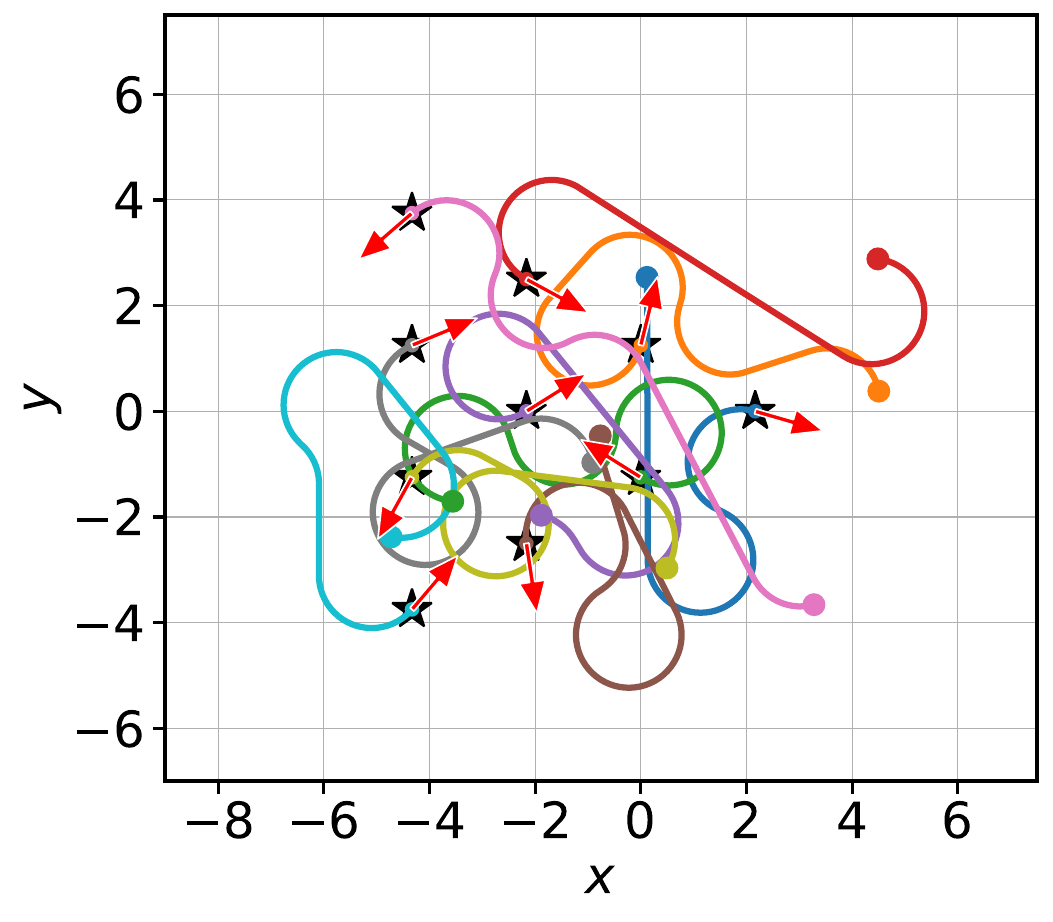}}\\
  \end{tabular}
  \caption{Formation planning with free terminal heading angles using the residual map implied by Theorem~\ref{thm:dubins2d}. The same instances and terminal times $T$ as in Figure~\ref{fig:Dubins2D_formation_fixed} are used. The red arrows indicate the obtained terminal directions.}
\label{fig:Dubins2D_formation_free}
\end{figure}
\subsection{Curvature-bounded paths in $\mathbb{R}^3$}
Similar to Dubins car, Theorem~\ref{thm:dubins3d} allows for applying nonlinear equation solvers to obtain the desired trajectory of given length. The difference is that no algebraic expression is available for H curves, and their endpoints must be calculated numerically. Thus, combination of \textsc{Matlab} \textit{ode45} and \textit{fsolve} functions is used to evaluate and find a zero of the error function that involves H curves. Throughout the simulations, the initial guesses were set by uniformly distributing the total time $T$ across the segments, and initializing the auxiliary parameters (e.g., orientation of the osculating plane and initial torsion) with a fixed set of constants. We note that testing feasibility of a given boundary condition is a reachability problem, which is distinct from trajectory generation. Given that describing the reachable sets of $f_{3D}$ remains an open problem~(\cite{bae2025reachability}), we utilize a feasible choice of terminal time $T$ for these demonstrations.

\begin{figure}[t]
    \centering
    \subfigure[Fixed $\boldsymbol{y}_f$ (CSCCSC)]{\includegraphics[width=0.48\columnwidth]{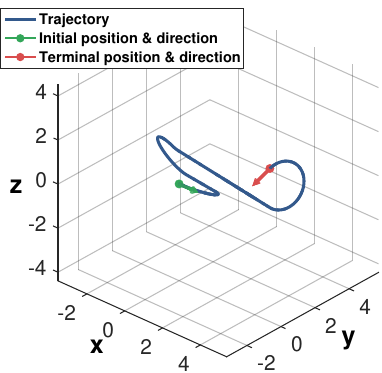}}
    \hfill
    \subfigure[Fixed $\boldsymbol{y}_f$ (HH)]{\includegraphics[width=0.48\columnwidth]{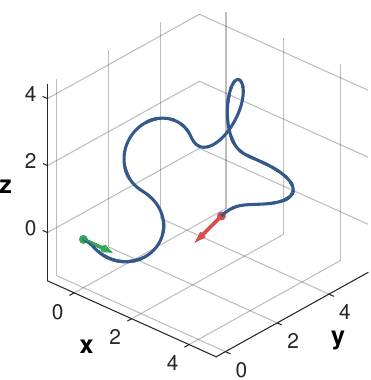}}
    
    \subfigure[Free $\boldsymbol{y}_f$ (CSCCSC)]{\includegraphics[width=0.48\columnwidth]{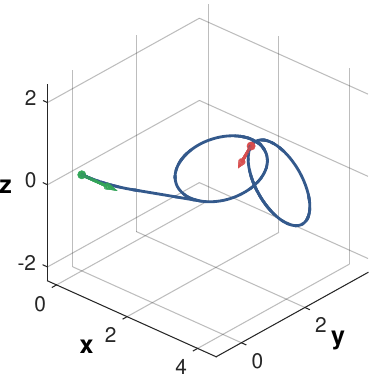}}
    \hfill
    \subfigure[Free $\boldsymbol{y}_f$ (HH)]{\includegraphics[width=0.48\columnwidth]{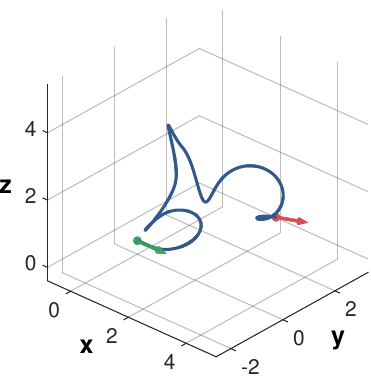}}
    
    \caption{Generated trajectories in $\mathbb{R}^3$ for fixed (top row) and free (bottom row) terminal velocity conditions. 
    Blue and red arrows indicate initial and terminal velocity vectors, respectively. 
    (a-b) Terminal velocity is fixed to $(\frac{1}{\sqrt{6}}, \frac{1}{\sqrt{3}}, \frac{1}{\sqrt{2}})$. 
    (c-d) Terminal velocity is unconstrained and is determined by the solver.}
    \label{fig:Dubins3D_combined}
\end{figure}

We examine the two scenarios of fixed and free terminal direction. In the latter case, the initial guess for the unconstrained terminal direction was determined based on the relative geometry between the boundary points. The specific boundary configurations and the obtained terminal velocities for the free case are summarized in Table~\ref{tab:3D_results}. The numerical procedure produces multiple solutions associated with different segment combinations characterized in Theorem~\ref{thm:dubins3d}. To illustrate this multiplicity, we present two representative solutions from different segment combinations, CSCCSC and HH, in Figure~\ref{fig:Dubins3D_combined}.
\begin{table}[tb]
    \centering
    \resizebox{\columnwidth}{!}{%
    \begin{tabular}{l c c c}
        \toprule
        \textbf{Inst.} & \textbf{Terminal Velocity} $\boldsymbol{y}_f$ & \textbf{Path Type} & \textbf{Resulting} $\boldsymbol{y}(T)$ \\
        \midrule
        \multirow{2}{*}{(a) Fixed $\boldsymbol{y}_f$} & \multirow{2}{*}{$\left(\frac{1}{\sqrt{6}}, \frac{1}{\sqrt{3}}, \frac{1}{\sqrt{2}}\right)$} & CSCCSC & Matches $\boldsymbol{y}_f$ \\
         & & HH & Matches $\boldsymbol{y}_f$ \\
        \midrule
        \multirow{2}{*}{(b) Free $\boldsymbol{y}_f$} & \multirow{2}{*}{\textit{Unconstrained}} & CSCCSC & $(0.44, -0.90, \phantom{-}0.02)$ \\
         & & HH & $(0.36, \phantom{-}0.85, -0.39)$ \\
        \bottomrule
    \end{tabular}
    }
	\caption{Simulation configurations and results for curvature-bounded paths in $\mathbb{R}^3$. Common parameters are set as follows: initial state $\boldsymbol{\chi}_{3D, i} = (0, 0, 0, 1, 0, 0)$, terminal position $\boldsymbol{x}_f = (3, 2, 1)$, and terminal time $T=6\pi$.}
    \label{tab:3D_results}
\end{table}

\begin{figure}[htbp]
	\centering
	\includegraphics[scale=1.0]{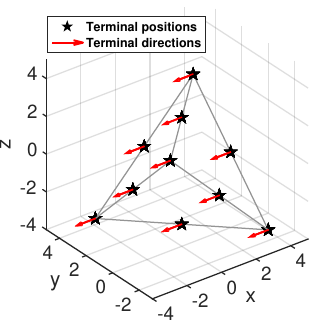}
	\caption{Setup of formation planning problem of curvature bounded paths in $\mathbb{R}^3$.}
	\label{fig:Dubins3D_formation_setup}
\end{figure}

\begin{figure}[htbp]
  \centering
  \setlength{\tabcolsep}{2pt}
  \begin{tabular}{@{}cc@{}}
    \scriptsize \textbf{Fixed $\boldsymbol{y}(T)$} & \scriptsize \textbf{Free $\boldsymbol{y}(T)$} \\[-0.1mm]

    \multicolumn{2}{@{}l@{}}{\scriptsize\textbf{Instance \#1}}\\[-0.5mm]
    \includegraphics[width=\imgw, trim=0 0 0 0]{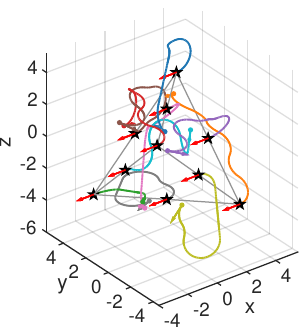} &
    \includegraphics[width=\imgw, trim=0 0 0 0]{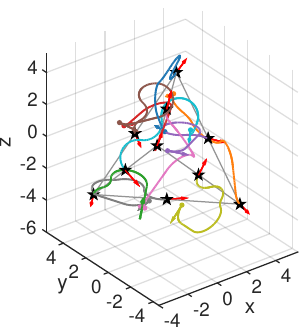} \\[-0.7mm]

    \multicolumn{2}{@{}l@{}}{\scriptsize\textbf{Instance \#2}}\\[-0.5mm]
    \includegraphics[width=\imgw, trim=0 0 0 0]{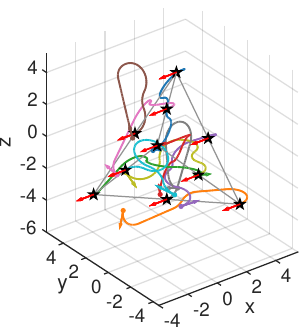} &
    \includegraphics[width=\imgw, trim=0 0 0 0]{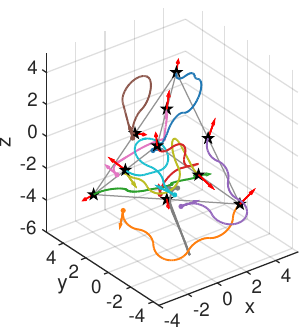} \\[-0.7mm]

    \multicolumn{2}{@{}l@{}}{\scriptsize\textbf{Instance \#3}}\\[-0.5mm]
    \includegraphics[width=\imgw, trim=0 0 0 0]{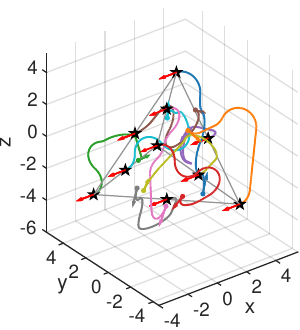} &
    \includegraphics[width=\imgw, trim=0 0 0 0]{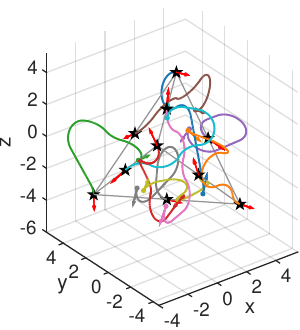} \\
  \end{tabular}
  \caption{Formation planning with fixed and free terminal heading angles for the three representative instances generated from random initial states. The red arrows indicate the terminal directions, which are prescribed (left) and obtained (right), respectively.}
  \label{fig:Dubins3D_formation}
\end{figure}

To further highlight the efficacy of Theorem~\ref{thm:dubins3d}, we address the 3D extension of the formation planning problem in Section~\ref{subsec:Dubins2D}. We conduct numerical experiments with 10 vehicles governed by the dynamics $f_{3D}$, tasked with forming the vertices and edge midpoints of a tetrahedron with an edge length of $8$ (see Figure~\ref{fig:Dubins3D_formation_setup}). Initial positions are drawn uniformly randomly from the set $[-3, 3]^3$, and initial directions from $S^2$. We visualize the results for three representative instances generated from these random initial states, considering both fixed and free terminal direction settings with a common terminal time of $T = 10$. Figure~\ref{fig:Dubins3D_formation} demonstrates that all vehicles successfully reach the prescribed terminal states, where the resulting terminal velocity vectors for the free scenario are visualized as red arrows. The exact numerical boundary values and these resulting terminal states are detailed in the provided code repository.
%
\section{Conclusion}\label{sec:06}
In this paper, we propose a formulation for the point-to-point steering problem of a broad class of nonlinear dynamics, which its Pontryagin extremals are computable. Theoretical completeness under the feasibility assumption is established: when the target state is reachable, we prove that a solution can always be constructed by concatenation of two Pontryagin extremals: one generated by the original dynamics $f$ from the initial state and the other generated by the inverted dynamics $-f$ from the terminal state. This yields a two-point boundary value problem (TPBVP) formulation whose solvability is equivalent to that of the original steering problem. The proposed formulation also accommodates cases where certain boundary states are left unconstrained. We further extend the framework to curves in $\mathbb{R}^3$ with a prescribed curvature bound, thereby extending the results on Dubins car to dimension three. In particular, it is proved that a solution to the associated steering problem can always be constructed by concatenating two curves of CSC, CCC, their subsegments, and H. Numerical demonstrations are provided for a controlled Van der Pol system, the Dubins car dynamics, and curvature-bounded paths in $\mathbb{R}^3$.
%
\begin{ack}
	The first author dedicates this work to the memory of his grandfather, Sang-Ki Park, whose unwavering support and encouragement have been a source of strength throughout his academic journey.

	This work was supported in part by the Institute for Information \& Communications Technology Planning \& Evaluation (IITP) grant funded by the Korea government (MSIT) (No. RS-2019-II190075, Artificial Intelligence Graduate School Program (KAIST)), and in part by the National Research Foundation of Korea (NRF) grant funded by the Korea government (MSIT) (No. RS-2026-25479034).
\end{ack}
\appendix
\section{Proof of Lemma~\ref{lemma:extendability}} \label{appendix:A}
\begin{pf*}{Proof.}
	The outline of the proof is as follows. We first show that the evolution of the state and costate variables of PMP (i.e., $(\boldsymbol{\chi}, \boldsymbol{p})$) can be equivalently described using a differential inclusion defined by the set-valued map $F(\boldsymbol{\chi},\boldsymbol{p})$. Then, we leverage the conclusions in~\cite{Filippov1988} for differential inclusions to establish the extendability of extremal trajectories, and compactness of $E(f,t_f,\boldsymbol{\chi}_i)$ follows from equicontinuity and Arzelà--Ascoli theorem.

	From continuity of $\mathcal{H}(\boldsymbol{\chi}, \boldsymbol{p}, \boldsymbol{u})$ in $\mathbb{R}^{2n_{\boldsymbol{\chi}}+n_{\boldsymbol{u}}}$ and compactness of $\Omega$, Berge's Maximum Theorem implies that $U(\boldsymbol{\chi}, \boldsymbol{p})$\textemdash{}and hence $F(\boldsymbol{\chi}, \boldsymbol{p})$\textemdash{}is compact and upper semicontinuous in $(\boldsymbol{\chi}, \boldsymbol{p})$. Let us consider the following differential inclusion.
	\begin{equation} \label{eq:lem2_diff_inc}
		\begin{bmatrix} \dot{\boldsymbol{\chi}} \\ \dot{\boldsymbol{p}} \end{bmatrix} \in F(\boldsymbol{\chi}, \boldsymbol{p})
	\end{equation}
	Suppose there is an extremal trajectory $\boldsymbol{\chi}(t)$ associated with costate trajectory $\boldsymbol{p}(t)$ and control input $\boldsymbol{u}(t) \in \mathcal{F}$. Then $(\boldsymbol{\chi}, \boldsymbol{p})(t)$ is absolutely continuous. Moreover, pointwise maximum Hamiltonian condition implies $\boldsymbol{u}(t) \in U(\boldsymbol{\chi}(t), \boldsymbol{p}(t))$ and hence $\begin{bmatrix} \dot{\boldsymbol{\chi}} \\ \dot{\boldsymbol{p}} \end{bmatrix} \in F(\boldsymbol{\chi}, \boldsymbol{p})$ almost everywhere. Thus, an extremal trajectory gives rise to a solution of the differential inclusion Eq.~\eqref{eq:lem2_diff_inc}.

	Conversely, suppose there is a solution $(\boldsymbol{\chi}, \boldsymbol{p})(t)$ of the differential inclusion Eq.~\eqref{eq:lem2_diff_inc}. If we write 
	\begin{equation}
		(\boldsymbol{\chi}, \boldsymbol{p})(t) = (\boldsymbol{\chi}, \boldsymbol{p})(0) + \int_0^{t} \phi(\tau) d\tau
	\end{equation}
	for an integrable function $\phi(t)$, then $\phi(t) \in F(\boldsymbol{\chi}(t), \boldsymbol{p}(t))$ for almost all $t$. Consequently, one can redefine $\phi(t)$ on a null set so that $\phi(t) \in F(\boldsymbol{\chi}(t), \boldsymbol{p}(t))$ for all $t$. From upper semicontinuity of $U(\boldsymbol{\chi}, \boldsymbol{p})$ and a gentle modification of Lemma 3A of Chapter~2 in~\cite{Lee1967}, there exists a measurable $\boldsymbol{u} \in \mathcal{F}$ such that $\phi(t) = \begin{bmatrix} f(\boldsymbol{\chi}(t), \boldsymbol{u}(t)) \\ -\nabla_{\boldsymbol{\chi}} f( \boldsymbol{\chi}(t), \boldsymbol{u}(t) )^T\boldsymbol{p}(t) \end{bmatrix}$ for all $t$. Hence, a solution of the differential inclusion in Eq.~\eqref{eq:lem2_diff_inc} gives rise to an extremal trajectory. Thus, the conclusions on differential inclusion theory in~\cite{Filippov1988} will be leveraged throughout the proof of this lemma.

	From the assumption, $F(\boldsymbol{\chi}, \boldsymbol{p})$ is compact, convex, and upper semicontinuous. Thus, $F(\boldsymbol{\chi}, \boldsymbol{p})$ satisfies the \textit{basic conditions} outlined in 2, $\S 7$, \cite{Filippov1988}. In the remaining steps, we will first show that the domain of $F(\boldsymbol{\chi}, \boldsymbol{p})$ can be restricted to an appropriate compact region. Then, conclusions in~\cite{Filippov1988} for differential inclusions on compact domains will be utilized.

	We can normalize the costate and further assume that $\|\boldsymbol{p}(0)\| = 1$ since the nontriviality condition is met by $\boldsymbol{p}$ itself. This restriction does not exclude any extremal trajectories since it is merely a normalization of the entire PMP conditions. The uniform bound assertion in Assumption~\ref{assum:01} implies the existence of a uniform bound as: $\left\| f(\boldsymbol{\chi}, \boldsymbol{u}) \right\| + \left\| -\nabla_{\boldsymbol{\chi}} f( \boldsymbol{\chi}, \boldsymbol{u})^T \right\| < \tilde{M}$. Consequently, it follows from costate differential equation and Grönwall's Inequality that 
	\begin{equation}
		\begin{split}
			\|\boldsymbol{p}(t)\| & \leq \|\boldsymbol{p}(0)\| + \int_0^t \left\| -\nabla_{\boldsymbol{\chi}} f( \boldsymbol{\chi}(\tau), \boldsymbol{u}(\tau))^T \right\| \|\boldsymbol{p}(\tau)\| d\tau \\
			& \leq 1 + \tilde{M} \int_0^t \|\boldsymbol{p}(\tau)\| d\tau \\
			& \leq e^{\tilde{M}t}.
		\end{split}
	\end{equation}
	Hence, $\boldsymbol{p}$ is uniformly bounded on $[0, 1]$, where the bound is independent of the choice of extremals. Thus, all extremals and the corresponding costate trajectories exist in an appropriate compact set $D$, and we can define $F(\boldsymbol{\chi}, \boldsymbol{p})$ on this domain to represent every extremal trajectories on $[0, T]$.

	If the basic conditions of $F$ are satisfied on a compact domain, any solution lying within the domain can be extended up to the boundary of the domain.~(Theorem~2, 2, $\S 7$, \cite{Filippov1988}) Therefore, any extremal trajectory can be extended on $[0, T]$. This proves the extendability claim.

	For compactness, let us consider a sequence of extremal trajectories on $[0, t_f] \subseteq [0, T]$ which their state, costate, and control are denoted as $\boldsymbol{\chi}_k(t)$, $\boldsymbol{p}_k(t)$, and $\boldsymbol{u}_k(t)$, respectively. Then each $(\boldsymbol{\chi}_k, \boldsymbol{p}_k)$ becomes a solution of the differential inclusion in Eq.~\eqref{eq:lem2_diff_inc} on $[0, t_f]$, lying in the compact region $D$. It then follows from Lemma~2 in 2, $\S 7$, \cite{Filippov1988} that that they are equicontinuous. Then the Arzelà-Ascoli theorem asserts the existence of a subsequence uniformly converging to a limit, $(\overline{\boldsymbol{\chi}}, \overline{\boldsymbol{p}})$. Since a uniform limit of solutions of the differential inclusion Eq.~\eqref{eq:lem2_diff_inc} is also a solution~(Corollary~1, 2, $\S 7$, \cite{Filippov1988}), it follows that $\overline{\boldsymbol{\chi}}$ is also an extremal. This completes the proof. \qed
\end{pf*}
%
\bibliographystyle{automatica}        
\bibliography{Ref_Papers}
\end{document}